\documentclass[11pt,a4paper]{article}
\usepackage{hyperref}
\usepackage[utf8]{inputenc}
\usepackage[T1]{fontenc}
\usepackage[margin=1in]{geometry}
\usepackage{amsmath,amsfonts,amssymb,amsthm}
\usepackage{physics,mathtools}
\usepackage{bbm}
\usepackage{cite}
\usepackage{xcolor,color,xspace}
\usepackage[noend]{algpseudocode}
\usepackage[ruled,vlined,linesnumbered]{algorithm2e}
\usepackage{caption, booktabs}
\usepackage{multirow}
\usepackage{url} 
\usepackage[capitalise]{cleveref}
\usepackage{array}

\newcommand{\R}{\mathbb{R}}
\newcommand{\Mr}{\mathcal{M}_r}
\newcommand{\Mri}{\mathcal{M}_{r_i}}

\newcommand{\cM}{\mathcal{M}}
\newcommand{\cB}{\mathcal{B}}
\newcommand{\cP}{\mathcal{P}}
\newcommand{\cO}{\mathcal{O}}
\newcommand{\cT}{\mathcal{T}}

\newcommand{\Matlab}{\textsc{Matlab}\xspace}
\newcommand{\Manopt}{\textsc{Manopt}\xspace}
\newcommand{\algcomment}[1]{\tcp{#1}}

\DeclareMathOperator{\diag}{diag}
\DeclareMathOperator{\vect}{vec}

\DeclareMathOperator{\sign}{sign}
\DeclareMathOperator*{\argmin}{arg\,min}

\newtheorem{theorem}{Theorem}[section]
\newtheorem{corollary}[theorem]{Corollary}
\newtheorem{lemma}[theorem]{Lemma}
\newtheorem{definition}[theorem]{Definition}

\newtheorem{remark}{Remark}

\title{Manifold-based Algorithms for the Hadamard 
Decomposition} 

\author{Nicolas Gillis\thanks{Department of Mathematics and Operational Research, University of Mons, Mons, Belgium. We acknowledge  the support by the European Union (ERC consolidator, eLinoR, no 101085607). \\ 
Emails: firstname.lastname@umons.ac.be.} \and Subhayan Saha\footnotemark[1] \and Stefano Sicilia\footnotemark[1] $^,$\thanks{Corresponding author. SS is a member of the Gruppo Nazionale Calcolo Scientifico-Istituto Nazionale di Alta Matematica (GNCS-INdAM).} \and Arnaud Vandaele\footnotemark[1]}

\begin{document}

 \maketitle
 
 \begin{abstract} 
  Given a matrix $X$, and two ranks $r_1$ and $r_2$, the Hadamard decomposition (HD) looks for two low-rank matrices, $X_1$ of rank $r_1$ and $X_2$ of rank $r_2$, both of the same size as $X$, such that $X\approx X_1\circ X_2$, where $\circ$ is the Hadamard (element-wise) product. 
  In most cases, HD is more expressive than 
  standard low-rank approximations such as the truncated singular value decomposition (TSVD), as it can 
  represent higher-rank matrices with the same number of parameters; this is because the rank of $X_1 \circ X_2$ is generically equal to $r_1 r_2$. 
  In this paper, we first present some theoretical insights for HD, in particular a useful reformulation $X\approx WH^\top$ where $W$ and $H$ have $r_1 r_2$ columns and belong to certain manifolds. These allow us to develop three new algorithms for computing HD. The first one uses the representation $X\approx X_1\circ X_2$ and relies on the Manopt toolbox. The other two rely on the reformulation $X\approx WH^\top$: one is a block projected gradient method, and the other is a manifold-based gradient descent algorithm that does not require projection onto the feasible set. The last two algorithms are particularly effective for handling large sparse data. We also propose new initializations that allow us to improve the accuracy of the HD. We compare our algorithms and initialization strategies with the TSVD and with the state of the art. Numerical results show that the new methods are efficient and competitive on both synthetic and real data. 
 \end{abstract}

  \vspace{0.1cm}
 

 \section{Introduction} 

 Finding a decomposition or factorization of a matrix $X$ is a well-known challenge in numerical analysis; see, e.g., \cite{tropp2017practical, tropp2019streaming, park2025low} for some recent papers on the topic. 
 This problem is generally motivated by the need to store the matrix $X$ efficiently, that is, to compress the dataset, especially when its dimensions are large. 
 Most of these models approximate the $m\times n$ matrix $X$ by using two matrices $W$ and $H$ of size $m\times r$ and $n \times r$, respectively, so that $X\approx WH^\top$. When $r\ll \min(m,n)$, this representation requires less memory to store the matrix $X$, and allows one to uncover some of its features and to provide computational advantages for subsequent data manipulation (e.g., fast matrix-vector product).  
 The most well-known model for low-rank compression is the truncated singular value decomposition (TSVD), essentially equivalent to principal component analysis (PCA). 
 It represents the original matrix $X$ as $X\approx U S V^\top$, where $U$ and $V$ have $r$ orthonormal columns and $S$ is a diagonal $r\times r$ matrix with nonnegative elements. The Eckart-Young theorem  ensures that the TSVD provides the best rank-$r$ approximation of $X$, both in norm 2 and in the Frobenius norm; 
 see, e.g., \cite{eckart1936approximation}. 

 Beyond data compression, the factors of the decomposition can identify important  features in the dataset $X$. In particular, if each column of $X$ represents a data point (e.g., a vectorized image or a document of word count), the columns of $W$ form a basis for the data points that can have a physical or probabilistic meaning. 
 Matrix factorizations have been used extensively in signal processing, data analysis, and machine learning; see, e.g., \cite{koren2009matrix, udell2016generalized, gillis2020nonnegative} and the references therein.

 Although many datasets are well-approximated by low-rank matrices~\cite{udell2019big, budzinskiy2025big}, 
 some may also have hidden features that are not easily captured by a linear model.  
 In order to overcome this limitation, other matrix decompositions have been considered, where the low-rank factors are combined using nonlinear functions. 
 In particular, Saul~\cite{saulnonlinear} recently introduced nonlinear matrix decompositions of the  
  form $X\approx f(WH^\top)$, where $f$ is a non-linear function applied componentwise to the entries of $WH^\top$. 
 Some examples are the following: $f$ is the ReLU function given by $f(x) = \max(0,x)$~\cite{saulnonlinear,seraghitiawarirelu,gillis2025extrapolated}, $f$ is the component-wise square function given by $f(x) = x^2$~ \cite{loconte2024turnknowledgegraphembeddings,csfjoakim},  $f$ is the sigmoid function given by $f(x) = \frac{1}{1+ e^{-x}}$~\cite{mnih2007probabilistic, sigmoidharrison}, 
 and $f$ is the min-max function, particularly well-suited for data within a bounded interval $[a,b]$, given by $f(x) = \min(b,\max(a,x))$~\cite{awari2025alternating}.

 In this work, 
 we focus on another type of nonlinear decomposition, namely the Hadamard decomposition (HD), that takes two low-rank matrices, $X_1$ and $X_2$, and combines them via the Hadamard  product, that is, the element-wise product denoted by $\circ$, so that $X_1 \circ X_2 \approx X$.  This decomposition mimics some real-world behavior, for instance in genetics (see, e.g.,  \cite{breslow1983multiplicative}) or in statistics, via connections to a family of graphical models called discrete restricted Boltzmann-machines (RBM); see \cite{oneto2023hadamard} for the connection between RBM and HD, and  \cite{fischer2012introduction} for a more general introduction to RBMs. 
 In signal processing, HD has been used to model the decomposition of an audio spectrogram $X$ into a low-rank spectral envelope part ($X_1$) and a low-rank pitch content ($X_2$) that interact through the element-wise product; see \cite{durrieu2011musically,durrieu2009main,durrieu2010source} for the details of this model called the   instantaneous mixture model. 
 In machine learning, HD has also been successfully used to compress or adapt neural networks~\cite{FedPara, huang2025hira}. 
 HD has also been extensively studied by the algebraic geometry community~\cite{onetoalggeom,kileelhad}. 
 In most cases, HD provides better compression than the TSVD, that is,  it provides a low-rank approximation with lower reconstruction error for the same number of parameters; see~\cite{ciaperoni2024hadamard, wertz2025efficient} and Section~\ref{sec:numexp} for numerical experiments. 

 Our focus in this paper is the practical computation of HD, which has not received much attention in the literature; only a few works have addressed this problem, namely  \cite{ciaperoni2024hadamard} and \cite{wertz2025efficient}. More specifically, we provide numerical optimization-based algorithms for computing an approximate solution to the HD problem: 
 given a matrix $X \in \R^{m \times n}$ and two ranks $r_1,r_2 \ll \min\{m,n\} $, solve
 \begin{equation} \label{eq:approxW1W2H1H2}   
  \min_{W_1,H_1,W_2,H_2} \left\|X - \left(W_1H_1^\top\right) \circ \left(W_2H_2^\top\right) \right\|_F^2, 
 \end{equation}
 where $W_1\in \R^{m\times r_1}$, $W_2\in \R^{m\times r_2}$, $H_1\in \R^{n\times r_1}$ and $H_2\in \R^{n\times r_2}$. 

\paragraph{Contribution and outline of the paper}

 The main contribution of this paper is to propose three new iterative algorithms for computing an HD, as well as effective initialization strategies. First, we provide some preliminary results on HD in Section~\ref{sec:theory} that are needed to develop the proposed approaches. 
 Section~\ref{sec:algorithms} presents our new algorithms for the computation of a HD based on two different manifold structures associated with the decomposition.  In Section~\ref{sec:Manopt}, we consider the standard structure, $X=X_1\circ X_2$, where $X_1$ and $X_2$ have fixed rank, and we rely on the \Manopt software to design an algorithm effective for small- and medium-sized matrices. In Section~\ref{sec:alternative}, we use an alternative structure that allows us to write the HD as a standard matrix factorization, $X=WH^\top$ under appropriate constraints on $W$ and $H$, yielding two approaches suited for large and sparse matrices:  Section~\ref{sec:projBCD} presents a simple two-block projected gradient descent scheme, while  Section~\ref{sec:manBCD} presents a more involved method that relies on a manifold-constrained gradient flow. 
 Section~\ref{sec:init} describes new initialization strategies to obtain good starting points for HD algorithms.  
 Section~\ref{sec:numexp} presents numerical experiments that compare our newly proposed algorithms and initialization strategies to the state of the art, showing the effectiveness of our novel methods on synthetic data, images, and real-world datasets. 

 \paragraph{Notation}

 For a matrix $A$, we use a \Matlab-like notation: we denote its $i$th row by $A(i,:)$, its $j$th column by $A(:,j)$ and its $(i,j)$ entry by $A(i,j)$, sometimes also indicated with the lowercase associated letter $a_{ij}$. The operator $\vect$ denotes the linear operation that converts an $m\times n$ matrix $A$ into the column vector $\vect(A)$ of length $mn$ obtained by stacking the columns of $A$ on top of one another; this operation is a bijection and its inverse is denoted as $\vect^{-1}$. 
 We use $\circ$ for the entry-wise multiplication, both for matrices and vectors, but also for the $k$th entry-wise power of a matrix $A$, which is denoted by $A^{\circ k}$. The transpose of $A$ is $A^\top$ and if $B$ has the same size as $A$,   their Frobenius inner product is $\langle A,B\rangle:=\trace(A^\top B)$, where $\trace(\cdot)$ is the trace, and the associated Frobenius norm $\|A\|_F^2=\langle A,A\rangle$. For any vector $v$, $\diag(v)$ is the diagonal matrix whose diagonal is $v$, and $\|v\|$ is the $\ell_2$-norm of $v$.

 \section{Preliminaries on the Hadamard decomposition} 
 \label{sec:theory}
 
 Let us introduce some definitions and present some known theoretical results concerning the Hadamard product of matrices. We rephrase and adapt these to our context in order to use them in the following sections. 

 \begin{definition}[Hadamard decomposition]
  \label{def:HD}
  Let $r_1$ and $r_2$ be two positive integers, and let $X\in~\R^{m\times n}$ be a given matrix. A rank-$(r_1,r_2)$ Hadamard decomposition (or factorization) of $X$ consists of a rank-$r_1$ matrix $X_1^\star\in \R^{m\times n}$ and a rank-$r_2$ matrix $X_2^\star\in \R^{m\times n}$ that solve
 \begin{equation}
  \label{eq:approxX1X2}
   \min_{(X_1,X_2)\in \cM_{r_1}\times\cM_{r_2}} \|X-X_1\circ X_2\|_F,  
 \end{equation}
 where $\circ$ denotes the Hadamard product, and $\Mri=\{M\in \R^{m\times n} : \rank(M)=r_i\}$.
 \end{definition} 
 We refer to $r_1$ and $r_2$ as the ranks of the Hadamard factors. 
 If $r_1=r_2=r$, we call the factorization a rank-$r$ HD. The HD is said to be exact if the minimum in \eqref{eq:approxX1X2} is $0$, that is, $X =X_1\circ X_2$. 
 Since the factor $X_i$ belongs to the rank-$r_i$ manifold $\Mri$, there exist factor matrices $W_i \in \R^{m\times r_i}$ and $H_i \in \R^{n\times r_i}$ such that $X_i=W_iH_i^\top$ for $i=1,2$. Hence \eqref{eq:approxX1X2} is equivalent to~\eqref{eq:approxW1W2H1H2}. The main goal of this paper is to provide novel efficient algorithms to solve this problem. 

 \paragraph{Alternative structure of the HD} 
 
Let us define the face-splitting product, which is closely related to the Khatri-Rao product. It will be useful for understanding the structure of the HD and to see it from a different perspective. 
 \begin{definition}[Face-splitting product]
  \label{def:facesplit}
  Let 
  \[
   W_1=\left(
   \begin{array}{c}
     a_1^\top \\
    \vdots \\
     a_m^\top
   \end{array}
   \right)\in\R^{m\times r_1}, \qquad W_2=\left(
   \begin{array}{c}
     b_1^\top \\
    \vdots \\
     b_m^\top
   \end{array}
   \right)\in \R^{m\times r_2}.
  \]
  Their face-splitting product is defined as 
  \[
  W_1\bullet W_2 := \left(
   \begin{array}{c}
     a_1^\top\otimes  b_1^\top \\
    \vdots \\
     a_m^\top\otimes  b_m^\top
   \end{array}
   \right)\in\R^{m\times (r_1 r_2)},
  \]
  where $\otimes$ denotes the Kronecker product.
  Equivalently, the face-splitting product can be defined using the columns of the matrices.
  Let 
  \[
   W_1=\left(
   \begin{array}{c|c|c}
     x_1 & \cdots &  x_{r_1}
   \end{array}
   \right)\in\R^{m\times r_1}, \qquad W_2=\left(
   \begin{array}{c|c|c}
     y_1 & \cdots &  y_{r_2}
   \end{array}
   \right)\in \R^{m\times r_2}, 
  \]
  then  
  \[
   W_1\bullet W_2 := \left(
   \begin{array}{c|c|c|c|c|c|c|c|c|c}
     z_{1,1} & \cdots &  z_{1,r_2} &  z_{2,1} & \cdots &  z_{2,r_2} & \cdots &  z_{r_1,1} & \cdots &  z_{r_1,r_2} 
   \end{array}
   \right)\in \R^{m\times (r_1 r_2)},
  \]
  where $ z_{i,j}= x_i\circ  y_j$ for all $i=1,\dots, r_1$ and $j=1,\dots, r_2$.
 \end{definition}
 
 It is straightforward to check that the two definitions of the face-splitting product coincide and that $(W_1\bullet W_2)^\top=W_1^\top \odot W_2^\top$, where $\odot$ denotes the Khatri-Rao product. Definition~\ref{def:facesplit} shows two different ways of building the matrix $W_1\bullet W_2$, highlighting its structure. The second definition states that each column of $W_1\bullet W_2$  is the entry-wise product of a column of $W_1$ and a column of $W_2$. The first definition shows that the $i$th row of $ W_1\bullet W_2$ is the vectorization of a rank-1 matrix, namely 
 \[
   (W_1\bullet W_2)(i,:) = \vect(b_i a_i^\top)^\top=a_i^\top \otimes b_i^\top=a_i^\top \bullet b_i^\top. 
 \]
This observation implies the following lemma.  
 \begin{lemma}
  \label{lem:W1fsW2}
  Let $W\in \R^{m\times (r_1 r_2)}$ and, for $k=1,\dots, m$, denote by $w_k^\top$ the $k$th row of $W$ and by $W^{(k)}$ the $r_2 \times r_1$ matrix such that $\vect(W^{(k)}) = w_k$. Then, there exist matrices $W_1\in \R^{m\times r_1}$ and $W_2\in \R^{m\times r_2}$ such that $W=W_1\bullet W_2$ if and only if $\rank(W^{(k)})\leq 1$ for all $k=1,\dots, m$.
  \begin{proof}
   It is a direct consequence of Definition~\ref{def:facesplit} of the face-splitting product.
  \end{proof}
 \end{lemma} 
 
 Now we recall a well-known result that describes how the face-splitting product is connected with the HD. 
 
 \begin{lemma}\cite[Lemma 3.1]{budzinskiy2025big},\cite[Lemma 4.16]{friedenberg2017minkowski},\cite[Proposition 1]{hyeon2021fedpara} 
  \label{lem:rankdec}
  Let $W_1\in \R^{m\times r_1}$, $W_2\in~\R^{m\times r_2}$, $H_1\in \R^{n\times r_1}$ and $H_2\in \R^{n\times r_2}$. 
  We have 
  \[
   (W_1H_1^\top)\circ (W_2H_2^\top)=\underbrace{(W_1\bullet W_2)}_{\in \R^{m\times (r_1 r_2)}}{\underbrace{(H_1\bullet H_2)}_{\in \R^{n\times (r_1 r_2)}}}^\top.
  \]
  Therefore, $\rank\left((W_1H_1^\top)\circ (W_2H_2^\top)\right)\leq r_1 r_2$.
  \begin{proof}
  First, we notice that given $ x, u\in \R^m$ and $ y, v\in \R^n$, it holds that
   \[
    ( x  y^\top)\circ( u  v^\top)=( x\circ  u)( y\circ  v)^\top,
   \]
   since the $(i,j)$ entry of both matrices is the product of the $i$th entries of $ x$ and $ u$ with the $j$th entries of $ y$ and $ v$.
   Now let $W_1(:,i)$, $H_1(:,i)$, $W_2(:,j)$, $H_2(:,j)$ denote the corresponding columns of $W_1,H_1,W_2$ and $H_2$, respectively.
   Then  
   \begin{align*}      
    (W_1H_1^\top)\circ (W_2H_2^\top) & = \left(\sum_{i=1}^{r_1}  W_1(:,i) H_1(:,i)^\top\right)\circ\left(\sum_{j=1}^{r_2}  W_2(:,j) H_2(:,j)^\top\right) \\ 
& 
    =\sum_{i=1}^{r_1} \sum_{j=1}^{r_2} \left(W_1(:,i) H_1(:,i)^\top\right)\circ \left(W_2(:,j) H_2(:,j)^\top\right) \\ 
    & =\sum_{i=1}^{r_1} \sum_{j=1}^{r_2} \left(W_1(:,i) \circ W_2(:,j)\right)\left(H_1(:,i) \circ H_2(:,j)\right)^\top=(W_1\bullet W_2)(H_1\bullet H_2)^\top, 
   \end{align*}
   where the last equality follows from the second definition of the face-splitting product.
  \end{proof}
 \end{lemma}

The following result further connects the Kronecker and face-splitting products to the HD. 
 \begin{lemma}\cite[Appendix A.3]{ciaperoni2024hadamard} 
  \label{lem:hadsvd}
  Let $X_1=U_1S_1V_1^\top$ and $X_2=U_2 S_2 V_2^\top$ be the thin TSVDs of $X_1,X_2\in \R^{m\times n}$ of rank $r_1$ and $r_2$, respectively. Then
  \[
   X_1\circ X_2=(U_1S_1V_1^\top)\circ (U_2 S_2 V_2^\top)=(U_1\bullet U_2)(S_1\otimes S_2)(V_1\bullet V_2)^\top.
  \]
  \begin{proof}
   For all $A\in \R^{m\times r_1}$, $B\in \R^{m\times r_2}$ and $M\in \R^{r_1\times r_1}$, $N\in \R^{r_2\times r_2}$ it holds that
   \[
    (A\bullet B)(M\otimes N)=\left(
   \begin{array}{c}
    (a_1^\top\otimes  b_1^\top)(M\otimes N) \\
    \vdots \\
    (a_m^\top\otimes b_m^\top)(M\otimes N)
   \end{array}
   \right)=\left(
   \begin{array}{c}
    (a_1^\top M)\otimes ( b_1^\top N) \\
    \vdots \\
    (a_m^\top M)\otimes ( b_m^\top N)
   \end{array}
   \right)=(A M)\bullet (B N),
   \]
   where $a_i^\top$ and $b_i^\top$ denote the rows of $A$ and $B$, respectively.
   Combining this identity with Lemma~\ref{lem:rankdec} gives
   \[
    (U_1S_1V_1^\top)\circ (U_2 S_2 V_2^\top)=((U_1 S_1)\bullet (U_2 S_2))(V_1 \bullet V_2)^\top=(U_1\bullet U_2)(S_1\otimes S_2)(V_1\bullet V_2)^\top.
   \]
  \end{proof}
 \end{lemma}

 \paragraph{Maximizing the expressivity of the HD}
 
 Lemma~\ref{lem:rankdec} implies $\rank(X)\leq r_1 r_2$ for an exact rank-$(r_1,r_2)$ HD of $X$. Moreover, such an HD has $(r_1+r_2)(m + n)$ parameters. To maximize the rank achievable by the HD, and hence the expressiveness of the model, it is therefore best to consider $r_1=r_2$. 
 The same observation was used in previous works~\cite{friedenberg2017minkowski, ciaperoni2024hadamard}. In the remainder of the paper, we assume $r_1 = r_2 = r$, unless stated otherwise. This also simplifies the presentation, although all results can be extended to the general case, $r_1 \neq r_2$. 
 
 One could also impose $X_1=X_2$ in HD~\eqref{eq:approxX1X2}, so that we get a decomposition represented by the componentwise square product, $X \approx X_1^{\circ 2}$. 
 This decomposition is useful, for instance, to design more expressive probabilistic circuits; see, e.g.,  \cite{awari2025alternating, loconte2022your}, or to compute the square 
 root rank~\cite{fawzi2015positive}.  
 However, this decomposition can only be used for a non-negative matrix $X$, and the rank of the resulting matrix, $X_1^{\circ 2}$, cannot reach $r^2$ for $r \geq 2$, as shown by the following corollary. 
 \begin{corollary}
  \label{lem:sqhad}
  Let $X=X_1\circ X_1$, where $X_1\in \R^{m\times n}$ has rank $r$. Then 
  \[
   \rank(X)\leq \frac{r(r+1)}{2}. 
  \]
  \begin{proof}
   Let $X_1=W_1H_1^\top$, with $W_1\in \R^{m\times r}$ and $H_1\in \R^{n\times r}$. Then, as done in Lemma~\ref{lem:rankdec}, we have
   \[
    (W_1H_1^\top)\circ (W_1H_1^\top)=\sum_{i=1}^r \sum_{j=1}^r (W_1(:,i)\circ W_1(:,j))(H_1(:,i)\circ H_1(:,j))^\top.
   \]
   In the summation, some rank-1 matrices share some vectors in common. The number of rank-1 matrices made up by distinct pairs of vectors is $\binom{r+1}{2}$, which yields the claim.
  \end{proof} 
   \end{corollary}

\paragraph{Existence and uniqueness of an  exact HD}

In general, it is not easy to decide if an exact HD exists. \cref{lem:W1fsW2} implies that there does not exist an exact rank-$r$ HD of a matrix of rank larger than $r^2$. 
If an exact rank-$r$ HD exists, then it is not unique. Indeed, one can always rescale the factors as follows 
 \[
  X= X_1 \circ X_2 \qquad X=(1/\alpha) X_1\circ \alpha X_2, \quad \forall \alpha\neq 0.
 \]
 Another way to generate multiple HDs is to leverage \cref{lem:W1fsW2}. Let $Z = uv^\top\in \R^{m\times n}$ be a rank-$1$ matrix with non-zero elements, that is, $u$ and $v$ have no zero entries. 
Let $WH^\top$ be a matrix of rank $r$ with  $W\in \R^{m\times r}$ and $H\in \R^{n\times r}$. Then
 \[
     (WH^\top) \circ Z = (WH^\top) \circ uv^\top = (W \bullet u) (H \bullet v)^\top = \tilde{W} {\tilde{H}}^\top, 
 \]
 where $\tilde{W} \in \R^{m\times r}$ and $\tilde{H} \in \R^{n\times r}$. Hence, the rank of $(WH^\top) \circ Z$ is not larger than $r$, and we always have the following non-uniqueness of HD: 
 \[
    X= X_1 \circ X_2 = (X_1 \circ Z) \circ (X_2 \circ Z^{\circ-1}),
 \] 
 where $Z^{\circ-1}$ is the component-wise inverse of $Z$.
 
 In this work, we are only interested in computing a HD; we leave the study of existence and uniqueness to future works. 

\paragraph{State-of-the-art algorithms for HD}

As far as we know, the study of practical algorithms for computing HD has received little attention in the literature. Here we briefly describe the two existing methods we are aware of. 

Ciaperoni et al.~\cite{ciaperoni2024hadamard} introduced the rank-$r$ HD in the context of data analysis and described a method that uses a gradient-based optimization approach to minimize 
\begin{equation}
    \label{eq:ciaperoni}
    E(W_1,H_1,W_2,H_2)=\|X - (W_1 H_1^\top)\circ(W_2 H_2^\top)\|_F^2.
\end{equation}
For a fixed learning rate $\eta$, they cyclically update each factor matrix (namely $W_1, H_1, W_2, H_2$) using a gradient descent step. 
To improve the convergence rate, they introduce a variant in which they multiply the gradients by scaling matrices. Let us show the update for the variable $W_1$, a similar procedure applies to the other factors. 
The scaling factor associated with the variable $W_1$ is $(H_1^\top H_1)^{-1}$, and its update is
\[
 W_1\leftarrow W_1-\eta \left(\nabla_{W_1} E\right) K, \qquad \nabla_{W_1} E=2\left(\left((W_1 H_1^\top)\circ(W_2 H_2^\top)-X\right)\circ (W_2 H_2^\top)\right)H_1,
\]
where $K=(H_1^\top H_1)^{-1}$ in the scaled version, and $K=I$ otherwise. The algorithm requires $\cO(mn\min(m,n)+Nmnr)$ operations, 
and $\cO(mn\min(m,n)+Nmnr+Nr^3+N\min(m,n)r^2)$ for the scaled variant, where $N$ is the number of iterations performed.

Wertz et al.~\cite{wertz2025efficient}  rely on  the same objective function \eqref{eq:ciaperoni} and use a block coordinate descent (BCD) method. Let us describe the update for the variable $H_2$,  the same approach applies for the other factors. 
For fixed $W_1,W_2$ and $H_1$, the subproblem with respect to $H_2$ can be written as   
\begin{equation*}
    \label{eq:wertz}
    \min_{H_2\in \R^{n\times r}} \sum_{j=1}^n \|X(:,j)-(W_1 H_1(j,:)^\top)\circ(W_2 H_2(j,:)^\top)\|_2^2, 
\end{equation*} 
showing that the rows of $H_2$ can be  independently updated.  Let us fix $j\in \{1,\dots,n\}$ and denote $x=H_2(j,:)^\top$. We need to solve the following least-squares problem:
\[
 \min_{x\in \R^r} \|(W_1 H_1(j,:)^\top)\circ (W_2 x) - X(:,j) \|_2^2, 
\]
whose associated normal equations are 
\[
  \left(W_2^\top \diag\left((W_1 H_1(j,:)^\top)^{\circ 2}\right) W_2\right) x=W_2^\top \left((W_1 H_1(j,:)^\top)\circ X(:,j)\right),
\]
which can be solved in $\cO(mr^2)$ operations. Hence, updating $H_2$ requires $\cO(m n r^2)$ operations, and similarly for the other factors. 
The authors accelerate their BCD algorithm using momentum, also known as extrapolation. We will use a similar strategy to accelerate our proposed algorithms in Section~\ref{sec:alternative}.  

These state-of-the-art methods do not exploit any manifold structure of the variables of the problem: our first main contribution is to describe a method that exploits the manifold pattern and is also able to use the information from the Hessian of the objective function to improve the HD provided (see Section~\ref{sec:Manopt}). Another drawback of the algorithms proposed in \cite{ciaperoni2024hadamard,wertz2025efficient} is that they cannot exploit the sparsity of $X$: indeed, they have to form $X_1$ and $X_2$ explicitly, which could be dense, possibly leading to memory issues. The other main contributions of this paper are to design two algorithms able to take advantage of the sparsity of $X$, allowing us to deal with larger matrices (see Section~\ref{sec:alternative}).

 \section{Manifold-based algorithms for the Hadamard decomposition}
 \label{sec:algorithms}

 In this section we propose three algorithms for computing a rank-$r$ HD of $X$. The first one relies on the standard representation, that is, 
 \begin{equation}
  \label{eq:HD1}
   X\approx X_1\circ X_2 =(W_1 H_1^\top)\circ (W_2 H_2^\top),  
 \end{equation}
 while the second and third ones exploit the alternative representation provided by Lemma~\ref{lem:rankdec}:
 \begin{equation}
  \label{eq:HD2}
   X\approx WH^\top =(W_1 \bullet W_2)(H_1\bullet H_2)^\top.  
 \end{equation} 
 In Section~\ref{sec:Manopt}, we follow the representation~\eqref{eq:HD1} to design a method that is well-suited for small- or medium-sized matrices; we refer to this algorithm as `Manopt', because it exploits the \Manopt software~\cite{boumal2014manopt}.  
 In Section~\ref{sec:alternative}, we use the alternative representation~\eqref{eq:HD2} to implement two scalable algorithms that are able to deal with large sparse matrices: in Section~\ref{sec:projBCD}, we propose a simple alternating projected gradient algorithm,  called `projBCD', while, in Section~\ref{sec:manBCD}, we propose a two-block coordinate descent method, called `manBCD', in which the variables are always kept on the manifold of the constraints.

 \subsection{Standard representation \texorpdfstring{$X=X_1\circ X_2$}{X=X1◦X2}: a second-order method for HD} 
 \label{sec:Manopt}
 
 We first consider the representation $X\approx X_1\circ X_2 =(W_1 H_1^\top)\circ (W_2 H_2^\top)$ in~\eqref{eq:HD1}. 
 The matrices $X_1$ and $X_2$ belong to the rank-$r$ manifold $\Mr$. This can be exploited by a Riemannian gradient descent (RGD) that looks for a pair $(X_1,X_2)\in \Mr\times \Mr$ that minimizes the functional
 \begin{equation}
  \label{eq:defphi}
  \Phi(X_1,X_2)=\frac{1}{2}\|X-X_1 \circ X_2 \|_F^2
 \end{equation}
 on the feasible set $\Mr\times \Mr$. To do so, we use the software \texttt{Manopt}~\cite{boumal2014manopt}, which, given an objective function, its Euclidean gradient (and possibly the Euclidean Hessian, which in this case can be computed efficiently; see below) 
 and a manifold for the constraints of the optimization, implements the RGD to compute a minimizer on the manifold. We need the expressions for the Euclidean gradient and the Euclidean Hessian of $\Phi$.
 \begin{lemma}
  \label{thm:phigradhess}
  Let $\Phi$ be defined as in~\eqref{eq:defphi} with the extended domain $\R^{m\times n}\times \R^{m\times n}$. Then, for all $X_1, X_2\in\R^{m\times n}$,
  \[
   \nabla \Phi(X_1,X_2)=
   \begin{pmatrix}
    -(X-X_1\circ X_2)\circ X_2 \\
    -(X-X_1\circ X_2)\circ X_1 \\
   \end{pmatrix}, 
  \]
  and 
  \[
   \nabla^2 \Phi(X_1,X_2)[A,B]=
   \begin{pmatrix}
    A\circ (X_2\circ X_2)+(2 X_1 \circ X_2-X)\circ B \\
    (2 X_1 \circ X_2-X)\circ A+B\circ (X_1\circ X_1) \\
   \end{pmatrix},
   \quad \forall A,B \in \R^{m\times n}.
  \]
  \begin{proof}
   Let us consider the vectorized versions of the matrices
   \[
    x=\vect(X); \quad u=\vect(X_1), \quad v=\vect(X_2),
   \]
   so that we can rewrite the objective function as
   \[
    \Phi(X_1,X_2):=\phi(u,v)=\frac{1}{2}\|x-u \circ v \|_2^2.
   \]
   Then, for $i=1,\dots, mn$, 
   \[
    \frac{\partial}{\partial u_i} \phi(u,v) = \frac{1}{2}\frac{\partial}{\partial u_i} \sum_{k=1}^{mn} (x_k-u_k v_k)^2=-(x_i-u_i v_i)v_i
   \]
   and, similarly, 
   \[
    \frac{\partial}{\partial v_i} \phi(u,v) = \frac{1}{2}\frac{\partial}{\partial v_i} \sum_{k=1}^{mn} (x_k-u_k v_k)^2=-(x_i-u_i v_i)u_i.
   \]
   Thus
   \[
    \nabla \phi(u,v)=-\begin{pmatrix}(x-u\circ v) \circ v\\(x-u\circ v) \circ u\end{pmatrix},
   \]
   which yields the first claim for the gradient of $\Phi$. For the computation of the Hessian of $\Phi$, we consider again the function $\phi$: for all $i,j=1,\dots, mn$, with $i\neq j$, it holds that
   \[
    \frac{\partial^2}{\partial u_i^2} \phi(u,v) = -\frac{\partial}{\partial u_i} (x_i-u_i v_i)v_i=v_i^2, \quad  
    \frac{\partial^2}{\partial u_i \partial u_j} \phi(u,v) = -\frac{\partial}{\partial u_i} (x_j-u_j v_j)v_j=0,
   \]
   \[
    \frac{\partial^2}{\partial v_i^2} \phi(u,v) = -\frac{\partial}{\partial v_i} (x_i-u_i v_i)u_i=u_i^2, \quad  
    \frac{\partial^2}{\partial v_i \partial v_j} \phi(u,v) = -\frac{\partial}{\partial v_i} (x_j-u_j v_j)u_j=0,
   \]
   \[
    \frac{\partial^2}{\partial u_i \partial v_i} \phi(u,v) = -\frac{\partial}{\partial u_i} (x_i-u_i v_i)u_i=2u_iv_i-x_i, \quad 
    \frac{\partial^2}{\partial u_i \partial v_j} \phi(u,v) = -\frac{\partial}{\partial u_i} (x_j-u_j v_j)u_j=0,
   \]
   and hence, since the function is smooth, we recover the remaining partial derivatives by means of Schwarz's theorem and get
   \[
    \nabla^2 \phi(u,v)=
    \begin{pmatrix}
     \diag(v\circ v) & \diag(2u\circ v-x) \\
     \diag(2u\circ v-x) & \diag(u\circ u) \\
    \end{pmatrix}.
   \]
   Finally, let $a=\vect(A)$ and $b=\vect(B)$ for $A,B\in \R^{m\times n}$. Then
   \[
    \nabla^2 \phi(u,v)
    \begin{pmatrix}
     a \\
     b \\
    \end{pmatrix}=
    \begin{pmatrix}
     \diag(v\circ v)a+\diag(2u\circ v-x)b & \diag(2u\circ v-x)a+\diag(u\circ u)b \\
    \end{pmatrix}^\top,
   \]
   which, after a proper reshaping and by recalling that $\diag(w)z=w\circ z$ for all vectors $w$ and $z$ of the same length, yields the expression for $\nabla^2 \Phi(X_1,X_2)[A,B]$.
  \end{proof}
 \end{lemma} 

 \paragraph{Implementation and computational cost} 

 Using the Euclidean gradient and Hessian of $\Phi$ from  Lemma~\ref{thm:phigradhess}, the implementation of an optimization algorithm on $\Mr\times \Mr$ is straightforward with \texttt{Manopt}, since the software computes the associated Riemannian gradient and Hessian. 

 The solver used in the implementation is \texttt{trustregions} and the overall cost of a single iteration of the method depends on the computation of the cost function $\Phi$, its gradient $\nabla \Phi$ and the action of the Hessian $\nabla^2 \Phi$. 
 As for the previous methods presented in Section~\ref{sec:theory}, Manopt requires forming $X_1$ and $X_2$ explicitly and therefore has to store $m\times n$ dense matrices while requiring $\cO(mnr)$ operations to build the matrices, which do not scale with the number of non-zeros in $X$. 
 For further details of the implementation, we refer to the code on Github  associated with the option `Manopt' in the main function `HadDec'; see  {\color{blue}\url{https://github.com/StefanoSicilia/Hadamard-Decomposition}}. 
 
 \subsection{Alternative representation \texorpdfstring{$X=WH^\top$}{X=WH'}: two scalable algorithms}
 \label{sec:alternative}
 
  We now propose two algorithms for HD based on the characterization of Lemma~\ref{lem:rankdec}:
 \begin{equation}
  \label{eq:decWH}
  X\approx WH^\top, \quad  W=W_1\bullet W_2, \quad H=H_1\bullet H_2.
 \end{equation} 
 In practice, on small- and medium-scale problems, these algorithms converge more slowly than Manopt proposed in Section~\ref{sec:Manopt}, but they can handle large sparse matrices, with a computational cost proportional to the number of non-zero entries in the input matrix.  

 Let us introduce the set 
 \begin{equation}\label{eq:defBmr}
     \cB_{m,r}=\left\{W=W_1\bullet W_2 : W_1, W_2\in \R^{m\times r}, 
     W(i,:) \neq 0 \text{ for all } i 
     \right\}\subseteq \R^{m\times r^2},
 \end{equation}
 whose closure is
 \[
  \bar{\cB}_{m,r}=\left\{W_1\bullet W_2 : W_1,W_2\in \R^{m\times r}  \right\}\subseteq \R^{m\times r^2}.
 \]
 The matrices $W$ and $H$ in~\eqref{eq:decWH} belong to $\bar{\cB}_{m,r}$ and $\bar{\cB}_{n,r}$, respectively. The  following three results provide some properties of these sets.
 \begin{lemma}
  \label{lem:Bmanifold}
  Let $m$ and $r$ be two positive integers such that $r\leq m$. Then the set $\cB_{m,r}$ is a manifold isomorphic to the product of $m$ copies of the vectorized rank-$1$ manifold $\cM_1\subseteq \R^{r\times r}$, 
  that is, 
  \[
   \cB_{m,r}\cong \prod_{i=1}^m \vect(\cM_1).
  \]
  \begin{proof}
   By Lemma~\ref{lem:W1fsW2}, if $B\in \cB_{m,r}$, then all its rows can be reshaped as rank-1 $r\times r$ matrices. This shows that $\cB_{m,r}$ is, up to a reshaping, the product of $m$ rank-1 manifolds $\cM_1$, one for each of its rows. 
  \end{proof}
 \end{lemma}

 \begin{lemma}
  Let $W, \widetilde{W}\in \bar{\cB}_{m,r}$. Then $W \circ \widetilde{W} \in \bar{\cB}_{m,r}$.
  \begin{proof}
    Let $x_iy_i^\top$ and $u_iv_i^\top$ be the rank-one matrices obtained by reshaping the $i$th row of $W$ and $\widetilde{W}$, respectively. 
    The corresponding row of $W\circ \widetilde{W}$ is $(x_iy_i^\top)\circ (u_iv_i^\top)=(x_i\circ u_i)(y_i\circ v_i)^\top$, which is still a matrix with rank not larger than $1$, meaning that $W\circ \widetilde{W} \in \bar{\cB}_{m,r}$. 
  \end{proof}
 \end{lemma}

 \begin{corollary}
  \label{cor:projBmr}
  Let 
  \[
   A=\begin{pmatrix}
    a_1^\top \\
    \vdots \\
    a_m^\top \\
   \end{pmatrix}\in \R^{m\times r^2}.
  \]
  The orthogonal projection of $A$ onto $\bar{\cB}_{m,r}$ with respect to the Frobenius norm is given by 
  \[
  \Pi_{\bar{\cB}_{m,r}}(A)=
   \begin{pmatrix}
    \vect\left({\Pi_{\cM_1}(\vect^{-1}(a_1))}\right)^\top \\
    \vdots \\
    \vect\left({\Pi_{\cM_1}(\vect^{-1}(a_m))}\right)^\top\\
   \end{pmatrix},
  \]
  where $\vect^{-1}$ is the $\vect$ inverse operator such that $\vect^{-1}(\vect(A))=A$, and $\Pi_{\cM_1}$ is the orthogonal projection onto $\cM_1\subseteq \R^{r\times r}$ given by the rank-$1$ TSVD. 
  \begin{proof}
   It is a direct consequence of Lemma~\ref{lem:Bmanifold} and the linearity of $\vect$.
  \end{proof}
 \end{corollary}

Algorithm~\ref{alg:projection} provides an implementation of the projection in Corollary~\ref{cor:projBmr}, exploiting the orthogonal projection onto rank-one matrices via the TSVD (Eckart-Young theorem).  

\begin{algorithm}[H]
\SetAlgoLined

\KwIn{An $m \times r^2$ matrix $A$}
\KwOut{The projection $W_1\bullet W_2$ of $A$ onto $\bar{\cB}_{m,r}$, defined in (\ref{eq:defBmr})
} 

\vspace{0.25cm}

Initialize $W_1$ and $W_2$ as empty matrices 

\For{$i = 1$ to $m$}{
    Reshape $A(i,:)$ as an $r \times r$ matrix and denote it by $A_i$ \\
    
    Let $(\sigma, u, v)$ be the rank-$1$ TSVD of $A_i$ \\
    
    Add $\sqrt{\sigma}u$ as a row to $W_1$ and $\sqrt{\sigma}v$ as a row to $W_2$ \\
}

\caption{$\Pi_{\bar{cB}_{m,r}}$, orthogonal projection onto $\bar{cB}_{m,r}$}
\label{alg:projection}
\end{algorithm}  
Algorithm~\ref{alg:projection} requires the computation of an SVD for each of the $m$ rows of the input matrix which are reshaped to form $r \times r$ matrices. Computing the SVD of an $r\times r$ matrix requires $\cO(r^3)$ arithmetic operations. Step~3 takes $\cO(r^2)$ arithmetic operations and computing the rows of $W_1$ and $W_2$ takes $\cO(r)$ arithmetic operations. 
The complexity of Algorithm \ref{alg:projection} is $\cO(m r^3)$. 

Before presenting our two proposed algorithms in the next two sections, let us provide the gradient and Hessian for the function we want to minimize, namely 
 \begin{equation}
  \label{prob:HadDecWH}
  \Psi : \R^{m\times r^2} \times \R^{n \times r^2} \to 
  [0,+\infty): (W,H) \mapsto 
  \Psi(W,H)=\frac{1}{2}\|X-WH^\top\|_F^2. 
 \end{equation}
 
 \begin{lemma}
  \label{lem:gradhesspsi}
  For the function $\Psi$ defined in~\eqref{prob:HadDecWH}, we have 
     \[
       \nabla \Psi(W,H)=
       \begin{pmatrix}
        \nabla_W \Psi(W,H) \\
        \nabla_H \Psi(W,H) \\
       \end{pmatrix}=
       \begin{pmatrix}
        (WH^\top-X)H \\
        (WH^\top-X)^\top W \\
       \end{pmatrix}:=
       \begin{pmatrix}
        G_W(W) \\
        G_H(H) \\
       \end{pmatrix} . 
      \]
     Moreover, for all $Y\in \R^{m\times r^2}$ and for all $Z\in \R^{n\times r^2}$, 
     \[
      \nabla^2_W \Psi(W,H) Y=Y H^\top H, \qquad \nabla^2_H \Psi(W,H) Z= Z W^\top W.
     \]
     \begin{proof}
         We prove the results for the block of variables $W$, the formulas for $H$ follow similarly. Let us consider the vectorized versions $x=\vect(X)$ and $w=\vect(W)$. Then, using the identities
         \[
          \|A\|_F^2=\|\vect(A)\|_F^2 \qquad \vect(AB)=(B^\top \otimes I_m)\vect(A), \quad \forall A\in \R^{m\times n}, \ \forall B\in \R^{n\times r},
         \]
         we have, for fixed $H$,
         \[
          \Psi(W,H) = \frac{1}{2}\|\vect(X)-\vect(WH^\top)\|_2^2=\frac{1}{2}\|(H\otimes I_m)w-x\|_2^2:=\psi(w).
         \]
         It is well-known that
         \[
          \nabla \psi(w)=(H\otimes I_m)^\top((H\otimes I_m)w-x), \qquad \nabla^2 \psi(w)=(H\otimes I_m)^\top (H\otimes I_m)=(H^\top H\otimes I_m).
         \]
 By reshaping, we get
         \[
          \nabla_W \Psi(W,H) = \vect^{-1}\left((H^\top\otimes I_m)(\vect(WH^\top-X))\right)=(WH^\top-X)H
         \]
         and 
         \[
          \nabla^2_W \Psi(W,H)Y=\vect^{-1}\left((H^\top H\otimes I_m)\vect(Y)\right)=Y H^\top H.
         \]
     \end{proof}
 \end{lemma}

The explicit expressions of $\nabla^2_W \Psi$ and $\nabla^2_H \Psi$ provide the Lipschitz constants of $\nabla_W \Psi$ and $\nabla_H \Psi$ (as univariate functions in $W$ and $H$ respectively), namely $L_W:=\|H^\top H\|_2$ and $L_H:=\|W^\top W\|_2$. These values will be used in the implementation of the algorithm to select a proper stepsize.

\subsubsection{Block Projected Gradient Descent} 
\label{sec:projBCD}

We first propose a two-block projected gradient descent scheme; see Algorithm~\ref{alg:projBCD}.  
The core idea is to alternatively update $W$ and $H$ by taking a step in the directions $-\nabla_W \Psi$ and $-\nabla_H \Psi$, and by projecting the result onto $\bar{\cB}_{m,r}$ and $\bar{\cB}_{n,r}$. 
We, however, adapt this simple idea in several ways to accelerate it: 
updating $W$ several times in a row, 
rescaling the factors at each iteration, and using extrapolation. We discuss the details of Algorithm~\ref{alg:projBCD} in the following paragraphs. We refer to the variable $W = W_1 \bullet W_2$ as the $W$-block, and similarly for $H$.

\begin{algorithm}[H]

\SetAlgoLined
\KwIn{A matrix $X \in \R^{m \times n}$, a positive integer $r$, a parameter $\tau$ for stepsize selection (default=0.95), a parameter $\beta_0$ as initial extrapolation parameter, number of inner iterations $k_W$ and $k_H$ (default=2)}
\KwOut{Matrices $W_1,W_2 \in \R^{m \times r}$, $H_1,H_2 \in \R^{n \times r}$ such that $X\approx(W_1H_1^\top)\circ(W_2H_2^\top)$}
\vspace{0.25cm}

Initialize $W_1 = W^{y}_1$, $H_1 = H^{y}_1$, $W_2 = W^{y}_2$, and $H_2 = H^{y}_2$  (see Section~\ref{sec:init}) \\   
Normalize $ X = \frac{X}{\|X\|_F}$ and initialize $\tilde{\beta} = 1, \beta = \beta_0$  \label{line:normalizeX} \\
 \While{some criterion is satisfied}{

\emph{\% Update of the $W$-block}

Let $(W_1^{y}, H^{p}_1, \text{vecnorm}(H_1^{y}))$ be the rescaling of $(W_1^{y}, H^{y}_1)$ following \eqref{eq:rescaling} \label{line:projbcd-W-start}\\

Let $(W_2^{y}, H^{p}_2, \text{vecnorm}(H_2^{y}))$ be the rescaling of $(W_2^{y}, H^{y}_2)$ following \eqref{eq:rescaling}\\ 

Set $H^{p} = H^{p}_1 \bullet H^{p}_2$,
$A = H^{p \top} H^{p}$, $B = X H^{p}$, $L_W = \|A\|_2$ and step-size $\alpha_W = \frac{\tau}{L_W}$ \label{line:precomputationW}

\For{$i = 1,2,\dots,k_W$}{ 
    Set $W^y = W^{y}_1 \bullet W^{y}_2$\\ 
    Update $G_W = W^y A-B$ \label{line:Wgradientcomputation}\\ 
    Update $(W^{y}_1,W^{y}_2)$ as output of Alg \ref{alg:projection} on $W^y-\alpha_W G_W$ \algcomment{Projection to $\bar{\cB}_{m,r}$} \label{line:updateW1W2}
}

Set $W_1^n = W_1^y \text{diag}({\text{vecnorm}(H_1^{y})}^{\circ -1})$, $W_2^n = W_2^y \text{diag}\left({\text{vecnorm}(H_2^{y})}^{\circ -1}\right)$ \algcomment{Scale} 

Set $W_1^y=W_1^n+ \beta (W_1^n-W_1)$, $W_2^y=W_2^n+ \beta (W_2^n-W_2)$ \algcomment{Extrapolation step} \label{line:projbcd-W-end}

\emph{\% Update of the $H$-block}  

Follow Lines~\ref{line:projbcd-W-start}--\ref{line:projbcd-W-end} to compute updates for $H_1^y$ and $H_2^y$ in exactly the same way. \\

    Set $H_1^y=H_1^n+ \beta (H_1^n-H_1)$, $H_2^y=H_2^n+ \beta (H_2^n-H_2)$ \algcomment{Extrapolation step}

    \emph{\% Checking decrease of the objective and adapting the  extrapolation parameter}  
    
    Set $W^n = W_1^n \bullet W_2^n$, $H^n = H_1^n \bullet H_2^n$  and compute the error $\|X - W^n H^{n \top}\|_F$  \label{line:errorcomputation} \\ 
    
   Update $\beta$ using Algorithm \ref{alg:betaupdate} \label{line:restartstart}\\
    \If{the error has increased}{ Set $W_1^y=W_1, W_2^y=W_2, H_1^y=H_1, H_2^y=H_2$ \algcomment{Extrapolation rejected} 
    \Else{Update $W_1=W_1^n, W_2=W_2^n, H_1=H_1^n, H_2=H_2^n$} }
\label{line:restartend}}
    \caption{Projection-based algorithm (projBCD) for computing a rank-$r$ HD}  
\label{alg:projBCD}
\end{algorithm}

\paragraph{Step-size Selection}

The function $\Psi$ is Lipschitz continuous in each block variable with Lipschitz constants 
\begin{equation}\label{eq:lipschitz}
    L_W = \|H^\top H\|_2 
    \quad \text{ and } \quad L_H = \|W^\top W\|_2.
\end{equation}
The step sizes $\alpha_W$ and $\alpha_H$ in the opposite direction of the gradient for the $W$ and $H$ block are chosen as $\alpha_W = \frac{\tau}{L_W}$ and $\alpha_H = \frac{\tau}{L_H}$, respectively, for some appropriate parameter $\tau \in (0, 2)$. We used $\tau=0.95$ in the numerical experiments.  
 
\paragraph{Multiple updates} 

In order to update the $W$-block, the matrices $A = H^\top H$, $B = XH$ and $L_W = \|A\|_2$ can be precomputed (see step \ref{line:precomputationW}). The second update of $W$ is therefore cheaper, and hence we propose to update the $W$-block $k_W$ times before updating the $H$-block. In the numerical experiments, we will use $k_W = k_H = 2$. 

\paragraph{Rescaling}\label{par:rescale}

Recall that for the $W$-block, the associated Lipschitz constant is $\| H^\top H \|_2$. To make the subproblems better conditioned, we resort to a standard rescaling procedure relying on the scaling degree of freedom in the product $WH^\top$. The procedure is applied in Step~\ref{line:projbcd-W-start} and does the following: given two matrices, $W_1^y \in \R^{m \times r}$  and $H_1^y \in \R^{n \times r}$, let 
\[
    \text{vecnorm}(H_1^y) := (\|H_1^y(:,1)\|_2,\dots,\|H_1^y(:,r)\|_2)^\top.
\]  
The procedure outputs the vector $\text{vecnorm}(H_1^y)$, the matrix $H_1^p$, and replaces the matrix $W_1^y$ where 
\begin{equation}\label{eq:rescaling}
    H_1^p = H_1^y \cdot \text{diag}(\text{vecnorm}(H_1^y))^{-1} 
    \quad \text{ and } \quad 
    W_1^y \leftarrow W_1^y \cdot \text{diag}(\text{vecnorm}(H_1^y)). 
\end{equation}
Essentially, it scales the $\ell_2$-norm of the columns of $H_1^y$ and $W_1^y$ so that $H_1^p$ has unit norm columns and $W_1^y(H_1^y)^\top$ is preserved. To avoid numerical issues, as a safety check, if a column of the matrix has a norm smaller than $ 10^{-15}$, we put $1$ in the corresponding entry of $\text{vecnorm}(H_1^y)$.

The rescaling procedure allows us to upper bound the Lipschitz constant in each block. Let $H^p_1$ and $H^p_2$ be the matrices obtained after the scaling procedure, that is, each of their columns has norm $1$. 
 
Since $H$ and $H^p$ are $n \times r^2$ matrices, we index their columns by $(i,j)$ where $i,j \in \{1,\dots,r\}$. Let $h$ and $h_p$ be the $(i,j)$-th column of matrices $H$ and $H^p$ respectively. Let $u$ and $u_p$ be the $i$-th column of $H_1$ and $H_1^p$  and let $v$ and $v_p$ be the $j$-th column of $H_2$ and $H_2^p$. We use the fact $\|h_{i,j}\|_2 = \|u \circ v\|_2 \leq \|u\|_2 \|v\|_2$ where the last inequality follows from the Cauchy-Schwarz inequality. We also have $u_p = \frac{u}{\|u\|_2}$, $v_p = \frac{v}{\|v\|_2}$, and $H^{p} = H^p_1 \bullet H^p_2$. Then the norm of $h^p_{i,j}$ can be bounded by  
\begin{align*}
    \left\|h^p_{i,j}\right\|_2 =\left\|u_p \circ v_p\right\|_2 &= \left\| \frac{u}{\|u\|_2} \circ  \frac{v}{\|v\|_2}\right\|_2  = \frac{\|h_{i,j}\|_2}{\|u\|_2 \|v\|_2} \leq 1.
\end{align*}
Using this, the Lipschitz constant after the scaling procedure is 
\[
    L_W = \|H^{p \top}H^{p}\|_2 \leq  \|H^{p \top}H^{p}\|_F \leq \sum_{i,j=1}^{r} \|h^p_{i,j}\|_2^2 \leq r^2.
\]

\paragraph{Extrapolation}\label{par:extrapolation}
The extrapolation strategy is inspired by the ones proposed in \cite{wertz2025efficient,ang2019accelerating}. 
Let $x$ denote our variable, and let  $x_{k+1}$ be the update for the $(k+1)$th step given by the gradient descent step on the $k$th iterate $x_k$. These iterates will tend to zigzag in  different directions, that is, $x_{k+1}-x_k$  and $x_k - x_{k-1}$ are usually very different directions. 
Adding extrapolation (a.k.a.\ momentum)  has a long history in optimization and allows one to obtain optimal first-order algorithms for convex problems; see~\cite{nesterov2013introductory} and the references therein. 
The idea is to introduce a new set of iterates, denoted by $y_k$ where $y_0 = x_0$ while $x_{k+1}$ is now the update given by the gradient step on $y_k$ and $y_{k+1} = x_{k+1} + \beta_k (x_{k+1} - x_k)$ where $\beta_k$ is the extrapolation parameter at the $k$th step.

Following this extrapolation strategy, for the update of the $W$-block, in Step~\ref{line:projbcd-W-end}, $W_1^n$ and $W_2^n$ are the updates given by the gradient step and the extrapolated updates are given by $W^y_1 = W_1^n  + \beta(W_1^n - W_1)$ and $W^y_2 = W_2^n  + \beta(W_2^n - W_2)$ where the variables $W_1$ and $W_2$ store the corresponding updates in the previous step. The extrapolated updates for the $H$-block are computed similarly.

\textit{Choice of the extrapolation parameter: } Algorithm~\ref{alg:betaupdate} implements a particular way to tune the extrapolation parameters $\beta$, inspired by the strategies in~\cite{wertz2025efficient,ang2019accelerating}. 

\begin{algorithm}[H]
\SetAlgoLined
\KwIn{The parameters $1 <  \tilde{\gamma} \leq \gamma \leq \eta, \beta_k, \beta_{\text{old}}$ and $\tilde{\beta}$}
\KwOut{The updated parameters $\beta_{k+1},\tilde{\beta}$ and $\beta_{\text{old}}$}
\vspace{0.25cm}
\If{the error has decreased in the current iteration}{set $\beta_{\text{old}} =     \beta$ \algcomment{Remember the last $\beta$ decreasing the error}
set $\beta_{k+1} = \min(\tilde{\beta},\gamma \beta_k)$ and $\tilde{\beta} = \tilde{\gamma}\tilde{\beta}$
}
\Else{Set $\tilde{\beta} = \beta_{\text{old}}$, $\beta_{k+1} = \frac{\beta_k}{\eta}$}
    \caption{Update of extrapolation parameters ($k$th step)}
    \label{alg:betaupdate}
\end{algorithm}

The algorithm starts with an initial value $\beta_0 \in [0,\tilde{\beta}]$, where the upper bound $\tilde{\beta} = 1$. When an iteration of Algorithm~\ref{alg:projBCD} decreases the error, the strategy is to increase the value of $\beta$ by a factor $\gamma > 1$, while taking into account the upper bound, that is, $\beta_{k+1} = \min(\tilde{\beta},\gamma\beta_k)$. Moreover, the upper bound is also increased by a factor $1 < \tilde{\gamma} < \gamma$. Let us consider the case where the error does not decrease with the current value of the extrapolation parameter denoted by $\beta$. In that case, we decrease the extrapolation parameter by the factor $\eta > \gamma$, and the upper bound $\tilde{\beta}$ is set to the previous value of $\beta$ that yielded a decrease of the error, that is, $\beta_{\text{old}}$.

In Steps \ref{line:restartstart}-\ref{line:restartend} in Algorithm \ref{alg:projBCD}, the current solution $(W_1,W_2,H_1,H_2)$ is updated only if the error has decreased; otherwise, the extrapolation step is rejected. In the rest of the paper, we will denote by $[\beta, \tilde{\beta}, \gamma, \tilde{\gamma}, \eta]$ to be the tuple of parameters required for the extrapolation step and the choice of the parameters for each of the algorithms will be further described in Section \ref{sec:effect}.

\paragraph{Complexity analysis of \texorpdfstring{Algorithm~\ref{alg:projBCD}}{Algorithm 2} }\label{par:companalprojBCD} 

Let $\sigma$ be the number of nonzero entries in $X$. In Step \ref{line:normalizeX} of the algorithm, the normalization of $X$ can be performed in time $\cO(\sigma)$. For the block of $W$ (Steps \ref{line:projbcd-W-start} to \ref{line:projbcd-W-end}),  Steps \ref{line:projbcd-W-start} and \ref{line:Wgradientcomputation} are the most expensive steps, which we explain in Table \ref{tab:complexity}. The total number of arithmetic operations required for each update of the $W$-block is $\cO(\sigma r^2+ (m+n)r^4 + r^6)$. Since the updates are symmetric in the $W$ and $H$ blocks, a similar analysis also holds for the update of the block of $H$.

The most expensive step in the rest of the algorithm is the computation of the error $\|X - W^n{H^n}^\top\|_F$ in Step~\ref{line:errorcomputation}. To do this, we exploit the fact that 
\[
 \|X - WH^\top\|_F^2 = \|X\|_F^2 - 2\langle X, W H^\top\rangle + \langle H^\top H,W^\top W\rangle. 
\]
For the second term, the $n \times r^2$ matrix $X^\top W$  can be computed in $\cO(\sigma r^2)$ arithmetic operations along with the fact that $\langle X, WH^\top\rangle = \langle X^\top W,H\rangle$. For the third summand, we need two $r^2 \times r^2$ matrices, $H^\top H$ and $W^\top W$, which can be computed in $\cO(nr^4)$ and  $\cO(mr^4)$ arithmetic operations respectively. 

\begin{table}[ht!]
\centering
    \resizebox{\textwidth}{!}{
    \begin{tabular}{|c|l|c|}
    \hline
    \textbf{Step} & \multicolumn{1}{c|}{\textbf{Operation}} & \textbf{Complexity} \\
    \hline
    
     & 
    Computation of $A = H^{p \top} H^{p} \in \mathbb{R}^{r^2 \times r^2}$ 
    & $\cO(nr^4)$ \\
    
    \ref{line:projbcd-W-start} & 
    Computation of $B = XH^{p}$ 
    & $\cO(\sigma r^2)$ \\
    
     & 
    Lipschitz constant $L_W = \|A\|_2$ (largest singular value of $A$) 
    & $\cO(r^6)$ \\
    
    \ref{line:Wgradientcomputation} & 
    Gradient $G_W = W^yA - B$ where 
    $W^y \in \mathbb{R}^{m \times r^2}$ and 
    $A \in \mathbb{R}^{r^2 \times r^2}$
    & $\cO(mr^4)$ \\
    \ref{line:errorcomputation} &  Computation of the error $\|X - W^n H^{n \top}\|_F$ & $\cO(\sigma r^2 + (m+n)r^4)$ \\
    \hline
    \end{tabular}
    }
    \caption{Complexity of Algorithm \ref{alg:projBCD} 
    }
    \label{tab:complexity}
\end{table}

Hence, the total number of arithmetic operations required for each iteration of Algorithm \ref{alg:projBCD} is $\cO(\sigma r^2 + (m+n)r^4 + r^6)$. Since typically $r \ll \min(m,n)$, the main cost comes from the term $\sigma r^2$, which scales linearly with the number of nonzero entries of $X$. If $X$ is dense, $\sigma = mn$ and this term becomes $mn r^2$.

\subsubsection{Gradient Flow on \texorpdfstring{$\cB_{m,r}$}{Bmr}}
 \label{sec:manBCD}

 Algorithm \ref{alg:projBCD} requires projection onto the manifold $\mathcal{B}_{m,r}$ as the gradient steps make the factor matrices leave the manifold. In this section, by analyzing a dynamical system from the gradient flows of a continuous version of the discrete variables $W$ and $H$, we propose a new algorithm whose updates remain on the manifold. The algorithm will also update $W$ and $H$ alternatively. It takes inspiration from a method developed for dynamical low-rank approximation \cite{koch2007dynamical} used for matrix nearness problems  \cite{guglielmi2025matrix,sicilia2025low}, and their applications in stabilization theory \cite{guglielmi2017matrix, guglielmi2023rank, guglielmi2024stabilization},  neural networks \cite{demarinis2025improving} and graph theory \cite{guglielmi2025low}. 
 The idea consists in reformulating the optimization problem as a gradient system whose flow is a continuous-time version of the variables to be optimized. The solution of the optimization problem corresponds to the stationary points of the gradient system, which turns out to be a matrix differential equation. 
 
 We adapt this approach to our setting. Let us introduce two matrix paths,  $W(t)$ and $H(t)$ with $t\geq 0$, that represent the two variables to be optimized. Then, taking into account the constraints, the gradient flows are given by  
 \begin{equation}
  \label{eq:gradsystWH}
  \dot{W}(t):=\dv{W(t)}{t}=-\cP_{W(t)} \left( G_W\left(W(t)\right) \right), \qquad \dot{H}(t):=\dv{H(t)}{t}=-\cP_{H(t)} \left(G_H\left(H(t)\right)\right),
 \end{equation}
 where $G_W$ and $G_H$ are the gradients of the objective function $\Psi$ with respect to $W$ and $H$, respectively (see Lemma \ref{lem:gradhesspsi}), $\cP_W$ is the projection onto the tangent space of $\cB_{m,r}$ at $W$, denoted as $\cT_W\cB_{m,r}$, and $\cP_H$ is the analogous projection for $\cT_H \cB_{n,r}$ at $H$; see  \cite{koch2007dynamical}. 
 Without the projections onto the tangent spaces, \eqref{eq:gradsystWH} consists of the free gradient system in each coordinate associated with the function $\Psi$. 

 Let us focus on the update of $W$ on the manifold $\cB_{m,r}$, the same approach applies to  $H$. 
 We first provide the expression of the projection $\cP_W$ onto the tangent space of  $\cB_{m,r}$ at $W$, relying on the simpler projection onto the tangent space to the manifold of rank-1 matrices.  
 \begin{lemma}
  \label{lem:tangprojBMr}
  Let $\cP_W$ be the orthogonal projection onto the tangent space $\cT_W\cB_{m,r}$ at 
  \[
   W=
   \begin{pmatrix}
    u_1^\top \otimes v_1^\top \\
    \vdots \\
    u_m^\top \otimes v_m^\top \\
   \end{pmatrix}=
   \begin{pmatrix}
    \vect(v_1u_1^\top) \\
    \vdots \\
    \vect(v_mu_m^\top)\\
   \end{pmatrix}=
   \begin{pmatrix}
    \vect(\rho_1 y_1x_1^\top) \\
    \vdots \\
    \vect(\rho_m y_mx_m^\top)\\
   \end{pmatrix}\in \cB_{m,r},
  \]
  where $\rho_i y_i x_i^\top=v_i u_i^\top$, $\rho_i\geq 0$ and $\|x_i\|=\|y_i\|=1$ for all $i=1,\dots,m$. For $A\in \R^{m\times r^2}$ with rows $a_i^\top$ for $i=1,\dots, m$, we have
  \[
   \cP_W(A)=
   \begin{pmatrix}
    \vect\left(A_1-(I_r-y_1 y_1^\top)A_1(I_r-x_1x_1^\top)\right)^\top \\
    \vdots \\
    \vect\left(A_m-(I_r-y_m y_m^\top)A_m(I_r-x_mx_m^\top)\right)^\top \\
   \end{pmatrix},
  \]
  where $A_i=\vect^{-1}(a_i^\top)$ and $I_r$ is the identity matrix of size $r$. Moreover, for all $A,B\in \R^{m\times r^2}$, 
  \begin{equation} \label{eq:projjorp} 
   \langle \cP_W(A), B \rangle = \langle \cP_W(A), \cP_W(B) \rangle =\langle A,\cP_W(B)\rangle. 
  \end{equation}
  \begin{proof}
   Let us consider the rank-1 manifold $\cM_1$ in $\R^{r\times r}$ and  denote by $\Pi_M$ the orthogonal projection onto the tangent space of $\cM_1$ at the point $M=\rho xy^\top$, where $x$ and $y$ have unit norm. We have  
   \[
    \Pi_M(B) = B - \left(I-xx^\top\right) B \left(I-yy^\top\right),
   \]
   which also means that $\langle A,\Pi_M(B)\rangle=\langle \Pi_M(A),B\rangle$ for all $A,B\in \R^{r\times r}$; see \cite{koch2007dynamical} or \cite[Proposition B.0.3]{sicilia2025low}. The expression of $\cP_W$ follows from Lemma~\ref{lem:Bmanifold},  reshaping the matrices into rows and by recalling that the tangent space of a product manifold is the product of the tangent spaces. Indeed every $W \in \mathcal{B}_{m,r}\cong \mathcal{B}_{1,r} \times \dots \times \mathcal{B}_{1,r}$ has a representation $(w_1,\dots,w_m) \in \mathcal{B}_{1,r} \times \dots \times \mathcal{B}_{1,r}$. 
   For the expression of $\cP_W$, we exploit  $\mathcal{T}_W\mathcal{B}_{m,r} \cong \times_{i=1}^m \mathcal{T}_{w_i}\mathcal{B}_{1,r}$.
   The equalities~\eqref{eq:projjorp} follow using Lemma~\ref{lem:Bmanifold} as well, since these equalities hold for $\Pi_M$.  
  \end{proof}
 \end{lemma}

  With the explicit expression for $\cP_W$, we propose an algorithm to compute the stationary points of~\eqref{eq:gradsystWH}, and discuss how these are related to the original optimization problem. 
  By omitting the dependence on $t$, we define $G_W:=G_W(W)$ and rewrite the first matrix ODE in~\eqref{eq:gradsystWH} as
 \begin{equation}
  \label{eq:odeW}
  \dot{W}=-\cP_{W} \left(G_W\right),
 \end{equation}
 where $\cP_{W} \left(G_W\right)$ is the Riemannian gradient of $\Psi$ with respect to $W$, and $G_W = \nabla_W \Psi$. Indeed, for a fixed $H\in \R^{n\times r^2}$, equation~\eqref{eq:odeW} is a gradient system for $\Psi_1(W):=\Psi(W,H)$, since along a solution $W(t)$, 
 \[
  \dv{t} \Psi_1(W(t))= \langle \dot{W}, G_W\rangle = -\langle \cP_W(G_W),G_W\rangle = -\|\cP_W(G_W)\|_F^2\leq 0.
 \]
 Thus, reaching a stationary point of system~\eqref{eq:gradsystWH} implies that the Riemannian gradient of the objective function is zero, that is,  we have found a stationary point of the function $\Psi_1$ on the manifold $\cB_{m,r}$. 

 We consider the problem row-wise, for $i=1,\dots,m$. For each row $ w_i^\top= u_i^\top\otimes  v_i^\top$ of $W$, we define its associated reshaped rank-$1$ matrix $M_i=v_i u_i^\top=\rho_i y_i x_i^\top$,  where $x_i = \frac{u_i}{\|u_i\|_2}$,  $y_i = \frac{v_i}{\|v_i\|_2}$,  and $\rho_i = \|u_i\|_2 \|v_i\|_2$. We rewrite the gradient system as
 \[
  \dot{M}_i=-\Pi_{M_i} G_i ,
 \] 
 where $G_i=G_i(x_i,y_i)=\vect^{-1}(G_W(i,:))\in \R^{r\times r}$ is the reshaped $i$th row of the gradient. By recalling the expression of $\Pi_{M_i}$ (see \cref{lem:tangprojBMr}), this is equivalent to
 \begin{equation}
    \label{eq:gradrhoxy} 
    \dot{\rho}_i y_i x_i^\top+\rho_i\dot{y}_i x_i^\top+\rho_i y_i\dot{x}_i^\top=-y_iy_i^\top G_i-G_ix_ix_i^\top+y_iy_i^\top G_ix_ix_i^\top.
 \end{equation}
 Equation~\eqref{eq:gradrhoxy} vanishes if $\rho_i$, $x_i$ and $y_i$ satisfy
 \begin{equation}
  \label{sys:gradrhoxy}
  \begin{dcases}
   \dot{\rho}_i=-(y_i^\top G_ix_i), \\
   \rho_i \dot{x}_i=-G_i^\top y_i+(y_i^\top G_ix_i) x_i, \\
   \rho_i \dot{y}_i=-G_ix_i+(y_i^\top G_ix_i) y_i. \\
  \end{dcases}
 \end{equation}
 This system also guarantees that the unit norm of $x_i$ and $y_i$ is preserved.
 Indeed
 \[
  \frac{1}{2}\dv{\|x_i(t)\|^2}{t}=x_i^\top \dot{x}_i=\frac{1}{\rho_i}\left(-(y_i^\top G_ix_i)^\top+(y_i^\top G_ix_i)\|x_i(t)\|^2\right)=0
 \]
 and similarly for $y_i$. Integrating~\eqref{sys:gradrhoxy} for all the rows would lead to the sought stationary points, but we prefer to integrate a slightly different system where the role of the scalar $\rho_i$ is replaced by the pair of norms $\mu_i = \|u_i\|$ and $\nu_i=\|v_i\|$. In this way it is possible to keep the balance between the norms of $u_i$ and $v_i$ even after the update, allowing the algorithm to converge to smaller relative errors in practice. Since $\rho_i=\mu_i \nu_i$, equation~\eqref{eq:gradrhoxy} is equivalent to
 \begin{equation}
    \label{eq:gradmunuxy} 
    \dot{\mu}_i \nu_i y_i x_i^\top+\mu_i \dot{\nu}_i y_i x_i^\top+\mu_i \nu_i\dot{y}_i x_i^\top+\mu_i \nu_i y_i\dot{x}_i^\top=-y_iy_i^\top G_i-G_ix_ix_i^\top+y_iy_i^\top G_ix_ix_i^\top,
 \end{equation}
 which vanishes if $\mu_i$, $\nu_i$, $x_i$, and $y_i$ satisfy
 \begin{equation}
  \label{sys:gradmunuxy}
  \begin{dcases}
   \dot{\mu}_i =-\lambda_i\frac{(y_i^\top G_ix_i)}{\nu_i}, \\
   \dot{\nu}_i =-(1-\lambda_i)\frac{(y_i^\top G_ix_i)}{\mu_i}, \\
   \mu_i \nu_i \dot{x}_i=-G_i^\top y_i+(y_i^\top G_ix_i) x_i, \\
   \mu_i \nu_i \dot{y}_i=-G_ix_i+(y_i^\top G_ix_i) y_i, \\
  \end{dcases}
 \end{equation}
 where $\lambda_i\in(0,1)$. In order to make the scheme symmetric, we select $\lambda_i=\frac{1}{2}$ for all $i$. Again, system~\eqref{sys:gradmunuxy} preserves the unit norms of $x_i$ and $y_i$.
 
 Now we briefly describe how to numerically integrate~\eqref{sys:gradmunuxy} when $\lambda_i=\frac{1}{2}$. Given the starting values $\mu_i^{(0)}$, $\nu_i^{(0)}$, $x_i^{(0)}$, and $y_i^{(0)}$, we use the explicit Euler method with a stepsize $h$ to update $x_i$ and $y_i$. Regarding the scalars $\mu_i$ and $\nu_i$, instead of the explicit Euler update, we use the exact solution of the first two equations where we consider $y_i^\top G_i x_i$ to be a constant. For $\mu_i$, we indeed have
 \begin{equation}
  \label{eq:muexact-euler}
   \mu_i(h)\approx\mu_i^{(0)}\sqrt{1-\frac{(y_i^\top G_ix_i)h}{\mu_i^{(0)}\nu_i^{(0)}}}=\mu_i^{(0)}-\frac{(y_i^\top G_ix_i)h}{2\nu_i^{(0)}}+\cO(h^2),  
 \end{equation}
 which results in a close approximation of the Euler step that, in practice, works better; an analogous equation also holds for $\nu_i$. Then, the integration scheme is
  \begin{equation}
  \label{sys:gradmunuxyeuler}
  \begin{dcases}
   \mu_i^{(k+1)}=\mu_i^{(k)} \omega_i^{(k)}, \qquad 
   \nu_i^{(k+1)}=\nu_i^{(k)}\omega_i^{(k)}, \\
   x_i^{(k+1)}=x_i^{(k)}+\frac{h}{\mu_i^{(k+1)}\nu_i^{(k+1)}} \left(-{G_i^{(k)}}^\top y_i^{(k)}+\left({y_i^{(k)}}^\top G_i^{(k)} x_i^{(k)}\right)x_i^{(k)} \right), \\
   y_i^{(k+1)}=y_i^{(k)}+\frac{h}{\mu_i^{(k+1)}\nu_i^{(k+1)}} \left(-G_i^{(k)} x_i^{(k)}+\left({y_i^{(k)}}^\top G_i^{(k)} x_i^{(k)}\right)y_i^{(k)} \right), \\
  \end{dcases} \qquad k=0,1,\dots
 \end{equation}
 where $\omega_i^{(k)}:=\sqrt{1-\frac{({y_i^{(k)}}^\top G_i^{(k)}x_i^{(k)})h}{\mu_i^{(k)}\nu_i^{(k)}}}$. During the integration, we select the same stepsize $h$ as for the projBCD method, inspired by the same Lipschitz constants $L_W=\|H^\top H\|_2$ associated with the objective function \eqref{prob:HadDecWH} of projBCD. In case $y_i^\top G_i x_i>0$, we also make sure that $h<\frac{\mu_i^{(0)}\nu_i^{(0)}}{(y_i^\top G_ix_i)}:=\bar{h}$ by selecting $h=\min(\frac{0.95}{L_W}, 0.95\bar{h})$ to avoid non-real solutions. This ensures that $\mu_i$ and $\nu_i$ are real and positive. Algorithm~\ref{alg:updmanifold} provides a single iteration of scheme~\eqref{sys:gradmunuxyeuler}. 
 
 The choice of a simple and explicit integration method is due to the fact that we are not interested in the exact trajectories of $x_i(t)$ and $y_i(t)$, but just in the stationary points towards which the dynamics converge. Moreover, since we update $x_i$ and $y_i$ directly, and consequently the components $u_i=\mu_i x_i$ and $v_i=\nu_i y_i$ of the rank-$1$ matrices corresponding to the $i$th row of $W$, by construction, we remain on the manifold $\cB_{m,r}$. Thus, we avoid the step of projecting back onto the manifold that is required by projBCD.

  \begin{algorithm}[ht!] 
  \SetAlgoLined
    \KwIn{Two matrices $W_1$ and $W_2$ of size $m\times r$, the gradient $G_W\in \R^{m\times r^2}$ and a stepsize $h$} 
    \KwOut{The updated $W_1$ and $W_2$}
\vspace{0.25cm} 
 \For{$i = 1$ to $m$}{
 Set $u_i=W_1(i,:)^\top$, $v_i=W_2(i,:)^\top$, $\mu_i=\|u_i\|$, $\nu_i=\|v_i\|$ and $\rho_i=\mu_i \nu_i$ \\   
    \If{$\rho_i\neq 0$}{
         Normalize to get $x_i=\frac{u_i}{\mu_i}$, $y_i=\frac{v_i}{\nu_i}$ \\
         Reshape $G_W(i,:)$ as an $r \times r$ matrix and denote it by $G_i$ \label{line:reshape}\\
         Precompute $g_i \xleftarrow[]{} G_i x_i$ and $\vartheta_i=y_i^\top g_i$ \label{line:gixi}\\
         If $\vartheta_i>0$ and $h>\frac{0.95}{\omega_i}$, set $h=\frac{0.95}{\omega_i}$. Then set $\omega_i=\sqrt{1-\frac{\vartheta_i h}{\rho_i}}$ \\
         Update $\mu_i \xleftarrow[]{} \mu_i \omega_i$, $\nu_i\xleftarrow[]{} \nu_i \omega_i$ and $\rho_i=\mu_i \nu_i$\\
         Update $x_i \xleftarrow[]{}  x_i+\frac{h}{\rho_i}\left(-G_i^\top y_i+\vartheta_i x_i\right) $  \label{line:xupdate}\\
         Update $y_i \xleftarrow[]{}  y_i+\frac{h}{\rho_i}\left(-g_i+\vartheta_i y_i\right) $   \\
         Normalize back by setting $W_1(i,:)=\mu_ix_i^\top$ and $W_2(i,:)=\nu_iy_i^\top$\\        
    } 
}
    \caption{Explicit step for system~\eqref{sys:gradmunuxy}} 
    \label{alg:updmanifold}
 \end{algorithm}

 The overall method, called `manBCD', results in a similar implementation as in Algorithm~\ref{alg:projBCD}, but Step~\ref{line:updateW1W2} is replaced by the following step
\begin{equation*}
    \text{Update } W_1^y,W_2^y \text{ as the output of Algorithm~\ref{alg:updmanifold} on input }( W_1^y,W_2^y, G_W, \alpha_W). 
\end{equation*} 
The update for $H_1^y,H_2^y$ is similarly replaced  by the output of Algorithm \ref{alg:updmanifold}.

 \paragraph{Complexity of \texorpdfstring{Algorithm~\ref{alg:updmanifold}}{Algorithm 4} and manBCD algorithm}

In Table \ref{tab:complexityupdman}, we extract the most expensive steps of Algorithm~\ref{alg:updmanifold} and provide their complexity analysis. 
\begin{table}[ht!]
\centering
\begin{tabular}{|c|c|c|}
\hline
\textbf{Step} & \centering\textbf{Operation} & \textbf{Complexity} \\
\hline

\ref{line:reshape} & 
Reshape $G_W(i,:)$ as an $r \times r$ matrix $G_i$
& $\cO(r^2)$ \\

\ref{line:gixi} & 
Computation of $g_i = G_ix_i$ 
& $\cO(r^2)$ \\

\ref{line:xupdate} & 
Computation of $G_i^\top y_i$
& $\cO(r^2)$ \\
\hline
\end{tabular}
\caption{Complexity of Algorithm \ref{alg:updmanifold}}
\label{tab:complexityupdman}
\end{table}
Hence, the total number of arithmetic operations for each row is $\cO(r^2)$, and the algorithm needs to perform this $m$ times. So the complexity of Algorithm~\ref{alg:updmanifold} is $\cO(mr^2)$.

From a complexity point of view, the same analysis as Algorithm \ref{alg:projBCD} follows, and one can get the same running time bounds, namely $\cO(\sigma r^2 + (m+n)r^4 + r^6)$ arithmetic operations, where $\sigma$ is the number of non-zero elements in the given matrix $X$.

 \paragraph{Discussion regarding the convergence guarantees} 
 
The restarting step of projBCD and manBCD ensures that the objective function is non-increasing at the end of each iteration of the algorithm. Since it is bounded below, the objective function values will converge. 

Regarding the convergence of iterates, Peng and Vidal~\cite{peng2023blockcoordinatedescentsmooth} provide a framework for BCD algorithms on manifolds. 
Algorithm~\ref{alg:projBCD} fits a framework analyzed in the paper, namely Rule (4): For the update of each block, take a step in the opposite direction of the gradient, then take a retraction of this step onto the manifold. But the corresponding analysis  \cite[Theorem 3]{peng2023blockcoordinatedescentsmooth} cannot be used in our case because it requires the corresponding manifold for each block to be a compact submanifold,  which does not hold in our case. 
Another recent work on the convergence of constrained block-Riemannian optimization methods is \cite{blockroptconv}. This could potentially be used to prove convergence guarantees for the iterates generated by Algorithm~\ref{alg:projBCD}. 
Unfortunately, their results also cannot be used in our context because one of the fundamental assumptions for the convergence guarantees, namely Assumption (A0) in \cite{blockroptconv}, is that the sub-level sets are compact which does not hold in our case.  We did not manage to find an existing framework or result that would prove the convergence of the iterates of our algorithms, and we leave this as a future direction of research.

 \subsection{Effect of acceleration and scaling} 
 \label{sec:effect}

Before reporting extensive numerical experiments in Section~\ref{sec:numexp}, we show the effect of acceleration and scaling on synthetic data. These two configurations cannot be applied to Manopt, and only the extrapolation is useful for BCD. Table~\ref{tab:comparison1} shows the relative error and the standard deviation (in percent) on 30 synthetic $400\times 400$ matrices\footnote{We use 10 matrices of each of the synthetic data matrices described in Section~\ref{sec:numexp}.} with $r=10,15,20$ and with a time limit of 10 seconds. 

\begin{table}[ht]
\centering
\begin{tabular}{c c c c c c}
\hline
$r$ & Method & none & extrapolation & rescaling & both \\
\hline
 & Manopt  & $\mathbf{44.78 \pm 0.06}$ & NA & NA & NA \\
 & manBCD  & $45.71 \pm 0.07$ & $45.85 \pm 0.08$ & $45.15 \pm 0.09$ & $\mathbf{44.89 \pm 0.08}$ \\
10 & projBCD & $46.02 \pm 0.07$ & $45.97 \pm 0.07$ & $45.42 \pm 0.08$ & $\mathbf{45.23 \pm 0.07}$ \\
 & BCD     & $44.95 \pm 0.06$ & $\mathbf{44.76 \pm 0.06}$ & NA & NA \\
 & TSVD    & $\mathbf{45.55 \pm 0.07}$ & NA & NA & NA \\
\hline
 & Manopt  & $\mathbf{42.13 \pm 0.05}$ & NA & NA & NA \\
 & manBCD  & $44.03 \pm 0.05$ & $43.90 \pm 0.13$ & $43.37 \pm 0.14$ & $\mathbf{42.58 \pm 0.11}$ \\
15 & projBCD & $44.62 \pm 0.04$ & $44.20 \pm 0.06$ & $43.65 \pm 0.10$ & $\mathbf{43.33 \pm 0.11}$ \\
 & BCD     & $42.47 \pm 0.06$ & $\mathbf{42.12 \pm 0.04}$ & NA & NA \\
 & TSVD    & $\mathbf{43.53 \pm 0.05}$ & NA & NA & NA \\
\hline
 & Manopt  & $\mathbf{39.44 \pm 0.07}$ & NA & NA & NA \\
 & manBCD  & $43.51 \pm 0.08$ & $42.56 \pm 0.08$ & $41.82 \pm 0.13$ & $\mathbf{40.95 \pm 0.16}$ \\
20 & projBCD & $44.19 \pm 0.09$ & $43.44 \pm 0.08$ & $42.05 \pm 0.08$ & $\mathbf{41.48 \pm 0.09}$ \\
 & BCD     & $39.96 \pm 0.07$ & $\mathbf{39.44 \pm 0.07}$ & NA & NA \\
 & TSVD    & $\mathbf{41.47 \pm 0.07}$ & NA & NA & NA \\
\hline
\end{tabular}
\caption{Effect of extrapolation and scaling on HD algorithms. 
The relative error and its standard deviation (in percent) are reported for each setting. 
For each configuration, a sample of $30$ synthetic matrices of size $400\times 400$ is generated.  with entries generated uniformly at random in [0,1]. An entry `NA' indicates that the corresponding configuration is Not Applicable. Best results per row in bold.}
\label{tab:comparison1}
\end{table}

As shown in \cite{wertz2025efficient}, the BCD method benefits from the extrapolation, as manBCD and projBCD do. Moreover, manBCD and projBCD improve their performance with the rescaling. In all cases, the best configuration for manBCD and projBCD is to use both extrapolation and rescaling. Hence, we will always set them both for manBCD and projBCD in the numerical experiments in Section~\ref{sec:numexp}, and we will always run BCD with extrapolation. Moreover, for all the experiments in the paper, the parameters for the extrapolation step, as described in paragraph~\hyperref[par:extrapolation]{\textbf{Extrapolation}} of Section~\ref{sec:projBCD}, are set to $[\beta, \tilde{\beta}, \gamma, \tilde{\gamma}, \eta] = [0.75,1,1.05,1.01,1.5]$ for BCD as recommended in~\cite{wertz2025efficient}, and to  $[0.25,1,1.05,1.01,1.5]$ for manBCD and projBCD (we also tried $\beta=0.5, 0.75$ which performed slightly worse than 0.25).

 \section{Initializations for the Hadamard decomposition}
 \label{sec:init}

Since computing a HD relies on iterative algorithms, choosing a good starting point may lead to better solutions and reduce the number of iterations. Only a few initialization strategies have already been proposed in the literature. In~\cite{ciaperoni2024hadamard}, the authors use random initialization. In~\cite{wertz2025efficient}, several initializations are proposed: random, SVD-based and $k$-means. The SVD-based initialization was reported to perform best. It decomposes $X$ as $X = X_1 \circ X_2$, where $X_1 = \sqrt{|X|}$ and $X_2 = \sign(X) \circ X_1$, and then replaces $X_1$ and $X_2$ with their rank-$r$ TSVD; see Algorithm~\ref{alg:initSVDbased}. 

 \begin{algorithm}[H]
 \SetAlgoLined
     \KwIn{The matrix $X$ and a positive integer $r$} 
     \KwOut{The matrices $W_1,W_2,H_1,H_2$ providing an initial approximation for the HD}

     \vspace{0.25cm}
     
      Compute $X_1=\sqrt{|X|}$ and $X_2=\sign(X) \circ X_1$, where $\sign(\cdot)$ and $\sqrt{\cdot}$ denote the entry-wise sign and square root function respectively. \\
      Let $(U_i,\Sigma_i, V_i)$ be the rank-$r$ truncated SVD of $X_i$ for $i=1,2$. \\
      Compute $W_i = U_i \sqrt{\Sigma_i}$ and $H_i = V_i \sqrt{\Sigma_i}$ for $i=1,2$. \\
    \caption{SVD-based Initialization for rank-$r$  HD~\cite{wertz2025efficient}} 
    \label{alg:initSVDbased}
 \end{algorithm}

We now propose three new SVD-based initializations. 
As we will see, they are in many cases better than the SVD-based initialization 
of~\cite{wertz2025efficient}. 
Their advantage is that they rely on the face-splitting structure highlighted in~\eqref{eq:decWH}: 
$X \approx W H^\top = (W_1 \bullet W_2) (H_1 \bullet H_2)$. 

The main idea behind the three proposed initialization strategies is to perform a rank-$r^2$ TSVD of $X$ to obtain $W$ and $H$, and then project $W$ and $H$ on their respective manifold, namely $\bar{\cB}_{m,r}$ and $\bar{\cB}_{n,r}$, to obtain a feasible solution, $(W_1,H_1,W_2,H_2)$.  

Hence, these initializations are well defined only if $r^2\leq \min(m,n)$, which is a condition that is generally satisfied in practice, since $r \ll \min(m,n)$ to obtain dimensionality reduction. 

We will discuss the practical performance of these strategies in the numerical experiments in Section~\ref{sec:numexp}.

 \paragraph{Initialization 1: Face-splitting (FS)}

As described in the previous paragraph, this first variant simply computes a rank-$r^2$ TSVD and projects the solution on the feasible set: 
 \begin{enumerate}
     \item Compute a TSVD of rank $r^2$ of $X$, $X \approx U \Sigma V^\top = \tilde U \tilde V^\top$, with $\tilde U = U \sqrt{\Sigma} \in \R^{m\times r^2}$ and $\tilde V = V \sqrt{\Sigma} \in \R^{n\times r^2}$.
     
     \item Project the factors $\tilde U$ and $\tilde V$ on the feasible sets: 
     \[
      W=W_1\bullet W_2=\Pi_{\bar{\cB}_{m,r}}(\tilde U), \qquad H=H_1\bullet H_2=\Pi_{\bar{\cB}_{n,r}}(\tilde V).
     \]
 \end{enumerate}

 \paragraph{Initialization 2 \& 3: 
 Face-Splitting-Left (FSL) \& Face-Splitting-Right (FSR) }
 
  Instead of simply taking $W$ and $H$ as the projection of the rank-$r^2$ TSVD, FSL only projects the second factor of the TSVD, $\tilde V$, while it recomputes the corresponding optimal left factor before its projection.  FSL performs the following 4 steps: 
  \begin{enumerate}
     \item Compute $(\tilde U, \tilde V)$ as in the first step of the FS initialization. 
     
     \item Project $\tilde V$ on the feasible set to obtain $H=H_1\bullet H_2=\Pi_{\bar{\cB}_{n,r}}(\tilde V)$.
     
     \item Compute $U_\star$ as the optimal solution of $\min_{U \in \R^{m\times r^2}} \|X- U H^\top\|_F$, that is, 
     $U_\star = X (H^{\dagger})^\top$, where $\dagger$ denotes the Moore-Penrose pseudo-inverse. 
         
     \item Project $U_\star$ on the feasible set to obtain $W=W_1\bullet W_2=\Pi_{\bar{\cB}_{m,r}}(U_\star)$. 
  \end{enumerate}

The FSR initialization is the FSL initialization applied to $X^\top$: first project $\tilde U$ to obtain $W$, then compute the optimal corresponding factor $V_\star$, and finally project it onto  $\bar{\cB}_{n,r}$ to obtain $H$.  

 \begin{remark}
   \label{rem:scaling_gamma}
     To improve any solution, such as the initializations discussed above, it is possible to rescale the factors, $W$ and $H$. 
     In fact, we can explicitly compute 
     \[
      \gamma_\star = \argmin_{\gamma\in \R} \|X-\gamma WH^\top\|_F^2 =\frac{\langle XH,W\rangle}{\langle W^\top W, H^\top H\rangle}, 
     \]
     and replace $(W,H)$ by $\sqrt{|\gamma_\star|}W$ and $\sign(\gamma_\star)\sqrt{|\gamma_\star|}H$.
 However, we did not observe significant improvements using this scaling on the initializations. Hence we did not include it in the initialization strategies. 
 \end{remark}

 \section{Numerical experiments}
 \label{sec:numexp}

 We now provide extensive numerical experiments comparing our proposed algorithms and initialization strategies with the state of the art. We compare the following algorithms: 
 \begin{itemize}
     \item Manopt described in Section~\ref{sec:Manopt} that relies on the \Manopt software. 

     \item projBCD described in Section~\ref{sec:projBCD} that uses a block projected gradient descent method. 
     
     \item manBCD described in Section~\ref{sec:manBCD}, which relies on gradient flows. 
     
     \item BCD: the block coordinate descent     
     method from \cite{wertz2025efficient}. The original implementation was in Julia, which we have translated into \Matlab and improved by removing unnecessary loops and by using vectorized operations that can be handled quickly by \Matlab. 
     
     \item TSVD: \Matlab \texttt{svd} to compute a TSVD matrix $X$. 
 \end{itemize}

\textit{Choice of extrapolation parameters:} 
As discussed in Section~\ref{sec:effect}, we use the following parameters from~\cite{wertz2025efficient} for BCD: $[\beta, \tilde{\beta}, \gamma, \tilde{\gamma}, \eta] = [0.75,1,1.05,1.01,1.5]$, and for projBCD and manBCD, we use $\beta = 0.25$ while keeping the other parameters the same. We do not include a comparison with the scaled gradient descent method proposed in \cite{ciaperoni2024hadamard}, since it is outperformed by BCD~\cite{wertz2025efficient}. 

Each time we run an algorithm for a rank-$r$ HD of $X$, we use our three initialization strategies presented in Section~\ref{sec:init}, along with the SVD-based initialization from \cite{wertz2025efficient}; see Algorithm~\ref{alg:initSVDbased}. 

\textit{Choice of measures of performance: }We report the lowest relative error, defined as $\frac{\|X-Y\|_F}{\|X\|_F}$ for a solution $Y$, among the four solutions obtained starting from these four initializations, and we will report which initialization achieved this lowest error. In case of a very small relative difference between the best initialization strategy and another one, namely less than $10^{-4}$, we consider it a tie and report both as best. 

To compare the compression performance of HD vs.\ TSVD, recall that a rank $2r$-TSVD uses the same number of parameters as a rank-$r$ HD. 
Let us denote $\text{err}_{\text{SVD}}^{r}(X)$ the relative error of the rank-$r$ TSVD of $X$, and 
$\text{err}_{\text{HD}}^{r}(X)$ the best found relative error of a rank-$r$ HD of $X$. 
We now define $r_\star$, which is intuitively the value such that the rank-$r_\star$ TSVD has a relative error closest to the rank-$r$ HD of $X$. 
More precisely, if $\text{err}_{\text{SVD}}^{2r}(X)<\text{err}^r_{\text{HD}}(X)$, we define 
\[
 r_\star = \max\left\{\varrho\in \mathbb{N} \ | \  \text{err}_{\text{SVD}}^\varrho(X)\geq \text{err}^r_{\text{HD}}(X)\right\}<2r,
\]
otherwise
\[
 r_\star = \min\left\{\varrho\in \mathbb{N} \ | \  \text{err}_{\text{SVD}}^\varrho(X)\leq \text{err}^r_{\text{HD}}(X)\right\}\geq 2r.
\] 
We can now define 
\[
 q_\star = \frac{r_\star - 2r}{2r},  
\] 
which measures the gain/loss of the compression ratio of HD compared to the TSVD. 
If $q_\star > 0$, it represents a gain of HD with respect to TSVD. 
For example, if $q_\star = 2$, it means that the TSVD needs twice as many parameters as HD to achieve the same relative error. 
If $q_\star < 0$, it is the other way around. As we will see, in most cases $q_\star > 0$. 

 All experiments have been performed on a laptop computer with Intel(R) Core(TM) i5-1035G1 CPU @ 1.00GHz (1.19 GHz) and 8GB of RAM with \Matlab R2024b. The code of the numerical experiments is available from  
 \begin{center}    
    {\color{blue}\url{https://github.com/StefanoSicilia/Hadamard-Decomposition}}. 
 \end{center}

We perform experiments on synthetic data (Section~\ref{sec:synt}), 
images (Section~\ref{sec:images}) and sparse document data sets (Section~\ref{sec:docs}).

 \subsection{Synthetic data} 
 \label{sec:synt}

 We analyze the performance of the compressions on randomly generated synthetic data of size $400\times 400$. We consider three types of randomly generated matrices: generic, low-rank and Hadamard-decomposable, using different values of the ranks, namely $r\in\{10,15,20\}$. 
 The generic random matrices are generated with entries drawn uniformly in [0,1] (using \texttt{rand(400)} in \Matlab). The low-rank type are rank-$2r$ matrices generated as the product of two matrices with inner dimension $2r$ whose entries are picked from the uniform distribution in [0,1], that is, \texttt{X = rand(400,2r) * rand(2r,400)} in \Matlab. 
For the Hadamard-decomposable type, we randomly generate matrices that admit a rank-$r$ HD: in \Matlab notation, 
  \[
   X = \left(\texttt{rand}(400,r) * \texttt{rand}(r,400) \right) 
  .* \left(\texttt{rand}(400,r) * \texttt{rand}(r,400)\right). 
  \]
 In each experiment, the number of  sampled matrices is fixed to 10. We run each algorithm with a time limit of 40 seconds for the generic case, and 100 seconds for the other two types. \cref{tab:trend} reports the results: the relative errors, along with the number of times each initialization performed best among the 10 samples (considering also ties). 

 \begin{table}[ht!]
     \centering
     \resizebox{\textwidth}{!}{
     \begin{tabular}{c c ccccc|c}
        \toprule
        & $r$ & Manopt & manBCD & projBCD & BCD & $2r$-TSVD
        & $\begin{array}{c} r_\star \\ q_\star (\%) \end{array}$ \\
        \midrule
        
        \multirow{6}{*}{\textbf{generic}}
        & \multirow{2}{*}{10} & $\underline{44.79 \pm 0.09}$ & $\underline{44.79 \pm 0.08}$ & $44.93 \pm 0.06$
        & $\mathbf{44.72 \pm 0.07}$ & $45.55 \pm 0.07$
        & \multirow{2}{*}{\begin{tabular}{c}$24.10 \pm 0.32$ \\ $20.50 \pm 1.58$ \end{tabular}} \\
        & & 0 - 10 - 0 - 0 & 0 - 10 - 0 - 0 & 6 - 5 - 0 - 0 & 6 - 4 - 0 - 0 &  & \\
        
        & \multirow{2}{*}{15} & $\underline{42.18 \pm 0.15}$ & $42.19 \pm 0.06$ & $42.71 \pm 0.04$
        & $\mathbf{42.04 \pm 0.05}$ & $43.53 \pm 0.05$
        & \multirow{2}{*}{\begin{tabular}{c}$38.00 \pm 0.00$ \\ $26.67 \pm 0.00$\end{tabular}} \\
        & & 0 - 10 - 0 - 0 & 0 - 10 - 0 - 0 & 5 - 5 - 0 - 0 & 6 - 4 - 0 - 0 &  & \\
        
        & \multirow{2}{*}{20} & $\underline{39.44 \pm 0.09}$ & $39.84 \pm 0.10$ & $40.86 \pm 0.10$
        & $\mathbf{39.28 \pm 0.07}$ & $41.47 \pm 0.07$
        & \multirow{2}{*}{\begin{tabular}{c}$52.00 \pm 0.00$ \\ $30.00 \pm 0.00$\end{tabular}} \\
        & & 0 - 10 - 0 - 0 & 0 - 10 - 0 - 0 & 0 - 10 - 0 - 0 & 10 - 0 - 0 - 0 &  & \\
        
        \midrule
        
        \multirow{6}{*}{\textbf{rank-$2r$}}
        & \multirow{2}{*}{10} & $\underline{1.21 \pm 0.04}$ & $\underline{1.21 \pm 0.04}$ & $\underline{1.21 \pm 0.04}$
        & $\underline{1.21 \pm 0.04}$ & $\mathbf{< 10^{-13}}$
        & \multirow{2}{*}{\begin{tabular}{c}$18 \pm 0$ \\ $-10 \pm 0$\end{tabular}} \\
        & & 0 - 10 - 0 - 0 & 0 - 10 - 0 - 0 & 0 - 10 - 0 - 0 & 5 - 6 - 0 - 0 &  & \\
        
        & \multirow{2}{*}{15} & $\underline{0.73 \pm 0.02}$ & $0.74 \pm 0.02$ & $0.74 \pm 0.02$
        & $0.74 \pm 0.02$ & $\mathbf{< 10^{-13}}$
        & \multirow{2}{*}{\begin{tabular}{c}$28 \pm 0$ \\ $-6.67 \pm 0$\end{tabular}} \\
        & & 0 - 10 - 0 - 0 & 0 - 10 - 0 - 0 & 1 - 9 - 0 - 0 & 0 - 10 - 0 - 0 &  & \\
        
        & \multirow{2}{*}{20} & $\underline{0.51 \pm 0.02}$ & $\underline{0.51 \pm 0.02}$ & $0.52 \pm 0.02$
        & $\underline{0.51 \pm 0.02}$ & $\mathbf{< 10^{-13}}$
        & \multirow{2}{*}{\begin{tabular}{c}$38 \pm 0$ \\ $-5 \pm 0$\end{tabular}} \\
        & & 0 - 10 - 0 - 0 & 0 - 10 - 0 - 0 & 0 - 10 - 0 - 0 & 0 - 10 - 0 - 0 &  & \\
        
        \midrule
        
        \multirow{6}{*}{\textbf{rank-$r$ HD}}
        & \multirow{2}{*}{10} & $< 10^{-5}$ & $0.08 \pm 0.24$ & $\mathbf{< 10^{-8}}$ & $\underline{< 10^{-7}}$ & $0.87 \pm 0.01$
        & \multirow{2}{*}{\begin{tabular}{c}$100 \pm 0$ \\ $400 \pm 0$\end{tabular}} \\
        & & 0 - 10 - 0 - 0 & 0 - 10 - 0 - 0 & 3 - 7 - 0 - 0 & 1 - 6 - 3 - 0 &  & \\
        
        & \multirow{2}{*}{15} & $\underline{< 10^{-6}}$ & $\mathbf{< 10^{-8}}$ & $0.41 \pm 0.09$
        & $< 10^{-2}$ & $0.60 \pm 0.01$
        & \multirow{2}{*}{\begin{tabular}{c}$225 \pm 0$ \\ $650 \pm 0$\end{tabular}} \\
        & & 0 - 10 - 0 - 0 & 0 - 10 - 0 - 0 & 2 - 8 - 0 - 0 & 1 - 9 - 0 - 0 &  & \\
        
        & \multirow{2}{*}{20} & $\mathbf{< 10^{-8}}$ & $0.28 \pm 0.11$ & $0.44 \pm 0$
        & $\underline{0.10 \pm 0.04}$ & $0.45 \pm 0.01$
        & \multirow{2}{*}{\begin{tabular}{c}$400 \pm 0$ \\ $900 \pm 0$\end{tabular}} \\
        & & 0 - 10 - 0 - 0 & 0 - 10 - 0 - 0 & 10 - 0 - 0 - 0 & 6 - 4 - 0 - 0 &  & \\
        
        \bottomrule
    \end{tabular}
    }
     \caption{Synthetic data: for each  of the nine experiments  (3 types of matrices: generic, rank-$2r$, and Hadamard decomposable, and 3 values for the rank $r \in \{10, 15, 20\}$), the first row reports the  relative errors in percent (average and standard deviation), and the second row reports number of times each initialization strategy is best (SVD-based - FS - FSL - FSR) with possible ties (up to $10^{-4}$ relative error). 
     It also reports the values of $r_\star$ and $q_\star$ (in percent). Best results are highlighted in bold, second best are underlined. 
     } 
     \label{tab:trend}
 \end{table}

 For the generic case, we observe that all the rank-$r$ HD outperform the $2r$-TSVD. For all the ranks, the best method is BCD, followed closely by Manopt, while the relative errors of manBCD and projBCD are slightly worse but comparable. The average relative gains of the HD ($q_\star)$ are respectively $20.5\%$, $26.67\%$, and $30\%$, showing a consistent advantage of the HD that grows as $r$ increases.

 In the low-rank case, the performances of the different methods for the HD are much closer, without significant differences. As expected, the $2r$-TSVD outperforms all rank-$r$ HDs, but the relative loss is limited, ranging from $-10\%$ to $-5\%$. The values of $r_\star$ are, on average, exactly $2r-2$, showing that the rank-$r$ HD can, in practice, approximate matrices of rank $2r$ well. 

 For the Hadamard-decomposable case, Manopt outperforms the other methods, being capable of recovering the exact decomposition for almost all the matrices. The other three methods are slightly worse, but their relative errors are better than those of the TSVD. The values of $r_\star$ match exactly the value of $r^2$, for which the low-rank decomposition always exists (see \cref{lem:rankdec}).

 Concerning the initialization strategy, the experiments show that it is not much related to the low-rank properties of the matrices, but more to the method: Manopt and manBCD always select FS.
 Besides FS, projBCD and BCD sometimes prefer the SVD-based initialization and more rarely FSL and FSR. 

 To summarize, BCD and Manopt are the best methods for synthetic datasets, with BCD being more accurate for the generic case, while Manopt is better for the Hadamard-decomposable matrices. The other two methods, manBCD and projBCD, perform a bit worse, but they still provide relative errors lower than the TSVD (with the obvious exception of the low-rank sample). The best initialization strategy is, on average, FS, followed by the SVD-based one.

\subsection{Dense images data sets} \label{sec:images}

\begin{figure}
    \centering
    \begin{minipage}{0.48\textwidth}
        \centering
        \includegraphics[width=\linewidth]{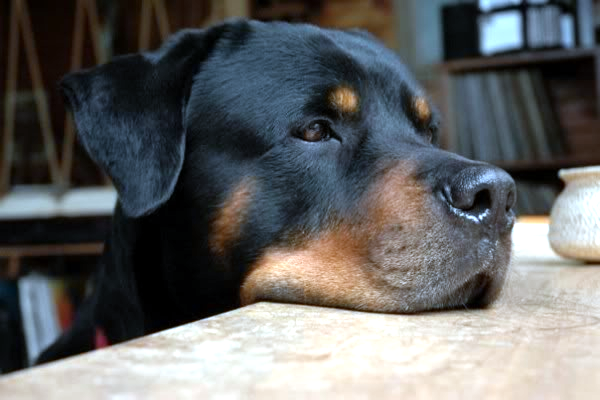}
        Original 
    \end{minipage}
    \hfill
    \begin{minipage}{0.48\textwidth}
        \centering
        \includegraphics[width=\linewidth]{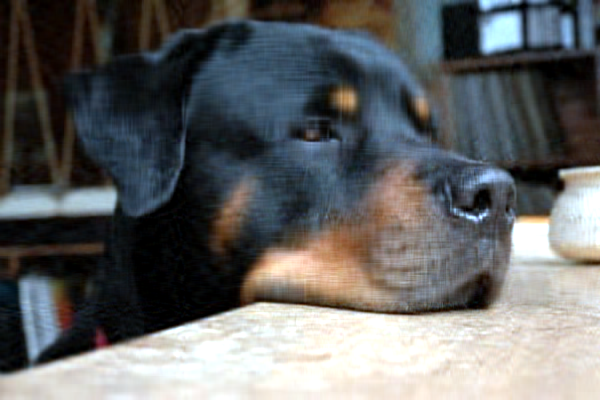}
       TSVD: 4.11\%, 4.85\%, 5.87\%
    \end{minipage}

    \vspace{0.3cm}
    
    \begin{minipage}{0.48\textwidth}
        \centering
        \includegraphics[width=\linewidth]{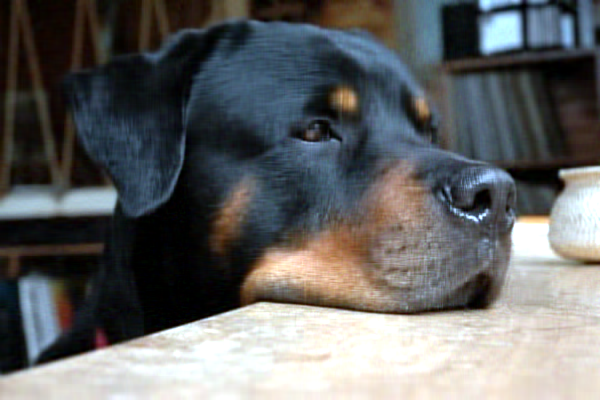}
        Manopt: 3.04\%, 3.61\%, 4.38\%
    \end{minipage}
    \hfill
    \begin{minipage}{0.48\textwidth}
        \centering
        \includegraphics[width=\linewidth]{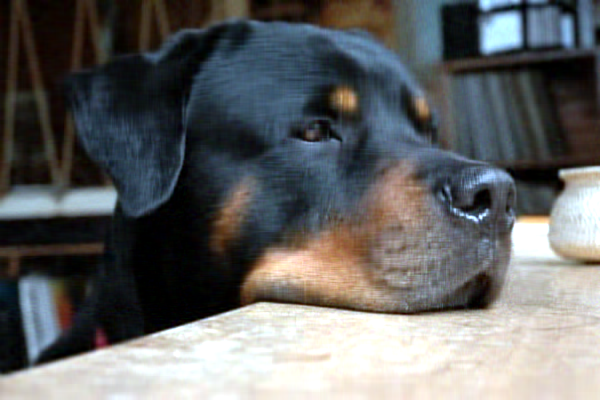}
        manBCD: 3.27\%, 3.88\%, 4.68\%
    \end{minipage}

    \vspace{0.3cm}

    \begin{minipage}{0.48\textwidth}
        \centering
        \includegraphics[width=\linewidth]{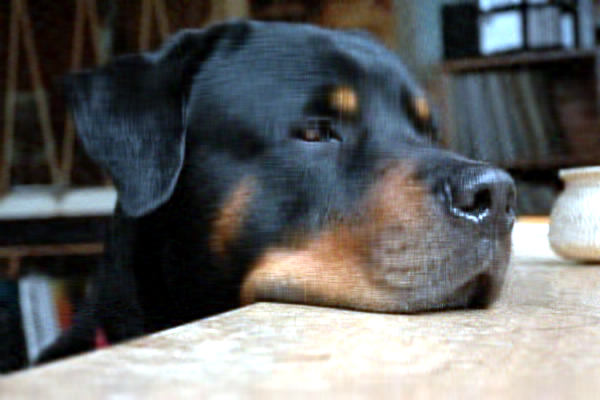}
        projBCD: 3.68\%, 4.31\%, 5.20\%
    \end{minipage}
    \hfill
    \begin{minipage}{0.48\textwidth}
        \centering
        \includegraphics[width=\linewidth]{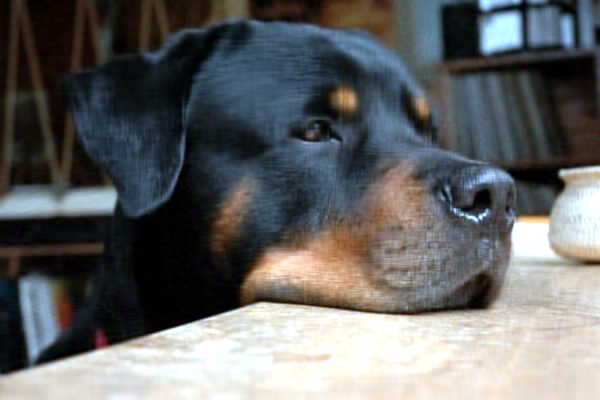}
        BCD: 3.18\%, 3.76\%, 4.57\% 
    \end{minipage}
    
    \caption{Dog image ($400\times 600\times 3$): reconstruction by  the different rank-20 HDs and by the rank-40 TSVD, relative error for each color band is reported. 
    In all cases the best initialization is FS, except for the second layer with BCD for which it is the SVD-based initialization.
    } 
    \label{fig:dog_approx}
\end{figure}

\begin{figure}[ht!]
\begin{tabular}{ccc}
    \includegraphics[width=0.3\linewidth]{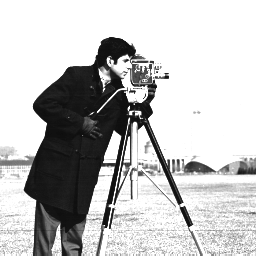} & \includegraphics[width=0.3\linewidth]{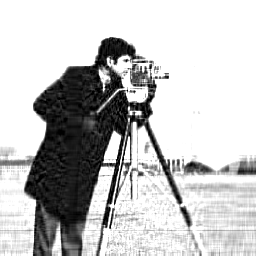} & \includegraphics[width=0.3\linewidth]{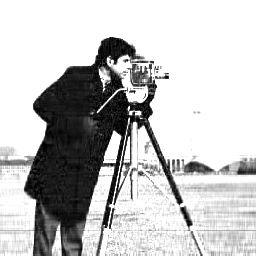}  \\
    Original & TSVD: 8.53\%  & Manopt: 6.06\% \\ 
    \includegraphics[width=0.3\linewidth]{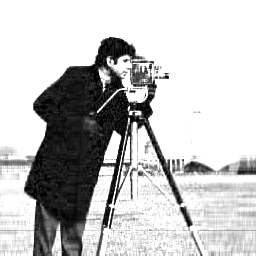} & \includegraphics[width=0.3\linewidth]{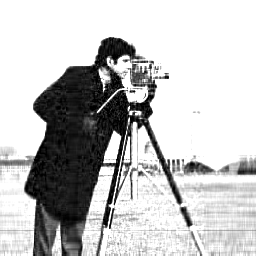} & \includegraphics[width=0.3\linewidth]{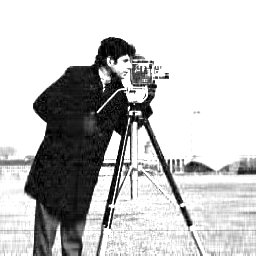} \\ 
    manBCD: 6.62\% & projBCD: 6.98\% &  BCD: 6.03\% 
\end{tabular}
    \caption{Approximations of the cameraman image ($256\times256$) with rank-16 HD and rank-32 TSVD. The relative error in percent is reported. 
    In all cases, the best initialization is FS. \label{fig:cameraman_approx} } 
 \end{figure}

We now perform experiments on images from \cite{ciaperoni2024hadamard}: the black-and-white Cameraman image (Figure~\ref{fig:cameraman_approx}), the cat image, three slices of the color dog image\footnote{In order to apply all the proposed initializations for $r=20$, we added an additional 400th row to the dog image which is originally of size $399 \times 600 \times 3$, copying the 399th row to have an image of dimension $400\times 600 \times 3$.}  (Figure~\ref{fig:dog_approx}), and the Olivetti faces datasets (a matrix containing vectorized $64\times 64$ images of faces). For each matrix, we consider the ranks $r \in \left\{\sqrt{\min(m,n)}/2, \sqrt{\min(m,n)}\right\}$. Table~\ref{tab:images} compares the different algorithms with a time limit of one minute per initialization.

\begin{table}[ht!]
 \centering
 \resizebox{\textwidth}{!}{
    \begin{tabular}{ccccccccccc}
        \hline
        Image & $m$ & $n$ & $r$ & Manopt & manBCD & projBCD & BCD & $2r$-TSVD & $r_\star$ & $q_\star$ (\%) \\
        \hline
        
        \multirow{2}{*}{cameraman} 
        & \multirow{2}{*}{256} & \multirow{2}{*}{256} 
        & 8  & \underline{10.78} (FS) & 11.02 (FS) & 11.20 (FS) & \textbf{10.62} (svd) & 12.98 & 23 & 43.75 \\
        & & & 16 & \underline{6.06} (FS) & 6.62 (FS) & 6.98 (FS) & \textbf{6.03} (FS) & 8.53 & 50 & 56.25 \\[0.5em]
        
        \multirow{2}{*}{cat} 
        & \multirow{2}{*}{400} & \multirow{2}{*}{600} 
        & 10 & \textbf{10.67} (FS) & 10.84 (FSL) & 11.05 (FSL) & \underline{10.76} (svd) & 12.74 & 32 & 60 \\
        & & & 20 & \underline{7.11} (FS) & 7.32 (FS) & 8.27 (svd) & \textbf{6.84} (FS) & 9.49 & 68 & 70 \\[0.5em]
        
        \multirow{2}{*}{dog (slice 1)} 
        & \multirow{2}{*}{400} & \multirow{2}{*}{600} 
        & 10 & \textbf{5.36} (FS) & 5.62 (FS) & 5.89 (FS) & \underline{5.56} (FS) & 7.07 & 29 & 45 \\
        & & & 20 & \textbf{3.04} (FS) & 3.27 (FS) & 3.68 (FS) & \underline{3.18} (FS) & 4.11 & 58 & 45 \\[0.5em]
        
        \multirow{2}{*}{dog (slice 2)} 
        & \multirow{2}{*}{400} & \multirow{2}{*}{600} 
        & 10 & \textbf{6.31} (FS) & 6.52 (FS) & 6.74 (FS) & \underline{6.44} (FS) & 8.20 & 29 & 45 \\
        & & & 20 & \textbf{3.61} (FS) & 3.88 (FS) & 4.31 (FS) & \underline{3.76} (svd) & 4.85 & 58 & 45 \\[0.5em]
        
        \multirow{2}{*}{dog (slice 3)} 
        & \multirow{2}{*}{400} & \multirow{2}{*}{600} 
        & 10 & \textbf{7.55} (FS) & 7.94 (FS) & 8.12 (FS) & \underline{7.57} (FS) & 9.66 & 29 & 45 \\
        & & & 20 & \textbf{4.38} (FS) & 4.68 (FS) & 5.20 (FS) & \underline{4.57} (FS) & 5.87 & 58 & 45 \\[0.5em]
        
        \multirow{2}{*}{olivettifaces} 
        & \multirow{2}{*}{4096} & \multirow{2}{*}{400} 
        & 10 & \textbf{11.55} (FS) & \underline{11.61} (FS) & 11.79 (FS) & 11.62 (svd) & 11.82 & 22 & 10 \\
        & & & 20 & \textbf{8.87} (FS) & \underline{9.11} (FS) & 9.48 (FS) & \underline{9.11} (svd) & 9.39 & 47 & 17.50 \\
        
        \hline
    \end{tabular}
    }
    \caption{Images: relative errors and best initialization for each algorithm. Best results are highlighted in bold, second best are underlined.} 
    \label{tab:images}
\end{table}

 In most of the cases, Manopt provides the lowest relative error. 
 The method is only outperformed by BCD for the cameraman image and the cat with $r=20$. Regarding the initialization strategy,  FS is most often the best, in more than $80\%$ of the cases. 
 The HDs outperform the TSVD in all cases, except projBCD for olivettifaces with $r=20$. 
 The relative gain of compression is significant, from the smallest $q_\star=10\%$ for the olivettifaces with $r=10$, to the highest 
 $q_\star=70\%$ for the cat image with $r=20$. These values highlight that, although the relative errors  do not differ much between HD and TSVD, there is a clear advantage in the compression rate, since the value of the rank $r_\star$ required by the TSVD to match the same accuracy of the HD is significantly larger than $2r$. This aspect is also confirmed by a clear visual improvement of the image reconstruction, especially in details. Figures~\ref{fig:dog_approx} and \ref{fig:cameraman_approx} provide such visual illustrations of the quality of the reconstructed images for the Cameraman image ($256\times 256$ pixels) for $r=16$ and the dog image for $r=20$, respectively. 
  For the dog image (Figure~\ref{fig:dog_approx}), we observe that all the rank-$20$ HDs perform better than the rank-$40$ TSVD, and the overall best approximation is provided by Manopt for all layers (relative errors of 3.03\%, 3.61\% and 4.38\%), which captures better details like the eyes and the nose of the dog. In all cases but one, our proposed FS initialization performs best. For all the slices, we find the same value of $r_\star=58$, which translates into a gain of $q_\star=45\%$ compared to the TSVD.  
 For the Cameraman image (Figure~\ref{fig:cameraman_approx}), Manopt and BCD provide very close approximations that are better than those of manBCD and projBCD (the man's coat is smoother, and the details of the camera are clearer). In any case, all HDs provide a much clearer reconstruction compared to the blurred  TSVD approximation. The initialization chosen is always FS, and the value of $r_\star$ is $50$, with a relative gain of the compression of $q_\star=56.25\%$.

\subsection{Sparse document datasets} \label{sec:docs} 

We now consider  text datasets from  \cite{zhong2005generative}, which are sparse matrices. The algorithms manBCD and projBCD exploit the sparsity structure, while Manopt and BCD cannot, as they both require building full $m\times n$ matrices during the computations (they need to explicitly form $X_1$ and $X_2$). All the algorithms are run for 200 seconds for each initialization, and Table~\ref{tab:text_data} reports the results. 

\begin{table}[ht!] 
    \centering

    \resizebox{\textwidth}{!}{
    \begin{tabular}{c c c c c c c c c c c}
        \hline
        Data & $m$ & $n$ & $r$ & Manopt & manBCD & projBCD & BCD & $2r$-TSVD & $r_\star$ & $q_\star$ (\%) \\
        \hline
        NG20    & 19949 & 43586 & 20 & -- & \textbf{9.50} (FS)  & \underline{10.48} (FSL) & -- & 10.29 & 71 & 77.5 \\
        classic &  7094 & 41681 &  4 & 93.77 (FSL) & \underline{90.81} (FS)  & \textbf{90.56} (svd) & -- & 92.68 & 14 & 75.0 \\
        la12    &  6279 & 31472 &  6 & 91.61 (FSR) & \textbf{86.19} (svd) & \underline{86.31} (svd) & -- & 89.78 & 26 & 116.7 \\
        ohscal  & 11162 & 11465 & 10 & 88.95 (FSR) & \textbf{81.14} (FS)  & \underline{81.33} (svd) & -- & 85.76 & 48 & 140.0 \\
        sports  &  8580 & 14870 &  7 & 86.32 (FS)  & \textbf{77.13} (FS)  & \underline{77.29} (svd) & 83.99 (FS) & 82.24 & 31 & 121.4 \\
        la1     &  3204 & 31472 &  6 & 88.40 (FSL) & \textbf{84.50} (svd) & \underline{84.78} (svd) & 88.02 (FSR) & 89.29 & 28 & 133.3 \\
        la2     &  3075 & 31472 &  6 & 87.82 (FSL) & \textbf{84.61} (svd) & \underline{84.95} (svd) & 87.99 (FSR) & 88.89 & 26 & 116.7 \\
        reviews &  4069 & 18483 &  5 & 81.49 (FS)  & \textbf{80.25} (FSL) & \underline{80.26} (svd) & 80.60 (FSL) & 83.30 & 20 & 100.0 \\
        k1b     &  2340 & 21819 &  6 & 84.32 (FSL) & \underline{83.48} (FS)  & 83.69 (svd) & \textbf{83.45} (svd) & 87.65 & 27 & 125.0 \\
        hitech  &  2301 & 10080 &  6 & 83.86 (FSR) & \underline{83.52} (FSL) & 83.66 (FSL) & \textbf{83.46} (svd) & 88.24 & 27 & 125.0 \\
        tr41    &   878 &  7453 & 10 & 31.05 (FS)  & \underline{29.93} (FS)  & 31.28 (svd) & \textbf{28.58} (svd) & 36.78 & 58 & 190.0 \\
        tr45    &   690 &  8261 & 10 & 14.02 (FS)  & \underline{13.27} (FSL) & 14.01 (svd) & \textbf{12.68} (svd) & 17.40 & 45 & 125.0 \\
        tr11    &   414 &  6424 &  9 & 13.89 (FS)  & \underline{13.62} (FS)  & 14.25 (svd) & \textbf{12.99} (svd) & 17.84 & 35 & 94.4 \\
        tr23    &   204 &  5831 &  6 & 10.14 (FS)  & \underline{9.91} (FS)   & 10.20 (svd) & \textbf{9.26} (svd) & 14.41 & 17 & 41.7 \\
        \hline
    \end{tabular}
    }
    \caption{Document data sets: relative errors and best initialization for each algorithm. An entry `--' denotes an error where \Matlab ran out of memory. Best results are highlighted in bold, second best are underlined. 
    } 
    \label{tab:text_data}
\end{table}

\newpage 

In all cases, HD offers a significant gain in compression compared to the TSVD (from 41.7\% to 190\%), most of the time with an impressive gain above 100\%. 
Manopt and BCD struggle with the largest datasets, and in some cases, they are not even able to provide a solution because they require too much memory. 
Note, however, that for the last six small datasets, BCD performs best. The experiment shows that for large sparse matrices, manBCD and projBCD are the methods of choice. Regarding the initializations, projBCD appears to perform best with the SVD-based initialization from~\cite{wertz2025efficient}, while the other algorithms more often perform best with one of our proposed schemes (FS, FSL, or FSR).

\section{Conclusion}
\label{sec:concl}

We proposed three new algorithms for computing rank-$(r_1,r_2)$ Hadamard decompositions (HD), relying on two different underlying manifold structures, as well as new initialization strategies. 
We have shown that these algorithms perform competitively with the state of the art and, in most cases, compress data matrices much better than the TSVD.  
The first method proposed, Manopt, which relies on the \Manopt software, 
performs best for dense medium-size matrices such as images,  
while the two other methods, manBCD and projBCD, 
are scalable for large sparse datasets such as documents. We also showed that the newly proposed initialization strategies are very successful in practice, often performing better than the state of the art, and that these initializations can also improve the performance of already existing methods, like BCD from~\cite{wertz2025efficient}. 

 {\small
 \bibliographystyle{abbrv}
 \bibliography{biblio}

@book{gillis2020nonnegative,
  title={Nonnegative matrix factorization},
  author={Gillis, Nicolas},
  year={2020},
  publisher={SIAM, Philadelphia}
}

@article{udell2019big,
  title={Why are big data matrices approximately low rank?},
  author={Udell, Madeleine and Townsend, Alex},
  journal={SIAM Journal on Mathematics of Data Science},
  volume={1},
  number={1},
  pages={144--160},
  year={2019},
  publisher={SIAM}
}

@inproceedings{huang2025hira,
  title={HiRA: Parameter-efficient {H}adamard high-rank adaptation for large language models},
  author={Huang, Qiushi and Ko, Tom and Zhuang, Zhan and Tang, Lilian and Zhang, Yu},
  booktitle={International Conference on Learning Representations},
  year={2025}
}

@article{FedPara,
  author       = {Na. Hyeon{-}Woo and M. Ye{-}Bin and T.{-}H. Oh},
  title        = {FedPara: Low-rank Hadamard Product Parameterization for Efficient Federated Learning},
  journal      = {ICLR}, 
  year         = {2022},
}

@article{mnih2007probabilistic,
  title={Probabilistic matrix factorization},
  author={Mnih, Andriy and Salakhutdinov, Russ R},
  journal={Advances in Neural Information Processing Systems},
  volume={20},
  year={2007}
}

@article{budzinskiy2025big,
  title={When big data actually are low-rank, or entrywise approximation of certain function-generated matrices},
  author={Budzinskiy, Stanislav},
  journal={SIAM Journal on Mathematics of Data Science},
  volume={7},
  number={3},
  pages={1098--1122},
  year={2025},
  publisher={SIAM}
}

@article{udell2016generalized,
  title={Generalized low rank models},
  author={Udell, Madeleine and Horn, Corinne and Zadeh, Reza and Boyd, Stephen and others},
  journal={Foundations and Trends{\textregistered} in Machine Learning},
  volume={9},
  number={1},
  pages={1--118},
  year={2016},
  publisher={Now Publishers, Inc.}
}

@article{tropp2017practical,
  title={Practical sketching algorithms for low-rank matrix approximation},
  author={Tropp, Joel A and Yurtsever, Alp and Udell, Madeleine and Cevher, Volkan},
  journal={SIAM Journal on Matrix Analysis and Applications},
  volume={38},
  number={4},
  pages={1454--1485},
  year={2017},
  publisher={SIAM}
}

@article{park2025low,
  title={Low-rank approximation of parameter-dependent matrices via CUR decomposition},
  author={Park, Taejun and Nakatsukasa, Yuji},
  journal={SIAM Journal on Scientific Computing},
  volume={47},
  number={3},
  pages={A1858--A1887},
  year={2025},
  publisher={SIAM}
}

@article{tropp2019streaming,
  title={Streaming low-rank matrix approximation with an application to scientific simulation},
  author={Tropp, Joel A and Yurtsever, Alp and Udell, Madeleine and Cevher, Volkan},
  journal={SIAM Journal on Scientific Computing},
  volume={41},
  number={4},
  pages={A2430--A2463},
  year={2019},
  publisher={SIAM}
}

@article{koren2009matrix,
  title={Matrix factorization techniques for recommender systems},
  author={Koren, Yehuda and Bell, Robert and Volinsky, Chris},
  journal={Computer},
  volume={42},
  number={8},
  pages={30--37},
  year={2009},
  publisher={IEEE}
}

@inproceedings{wertz2025efficient,
  title={Efficient algorithms for the hadamard decomposition},
  author={Wertz, Samuel and Vandaele, Arnaud and Gillis, Nicolas},
  booktitle={International Workshop on Machine Learning for Signal Processing (MLSP)},
  year={2025} 
}

@article{gillis2025extrapolated,
  title={An extrapolated and provably convergent algorithm for nonlinear matrix decomposition with the {R}e{LU} function},
  author={Gillis, Nicolas and Porcelli, Margherita and Seraghiti, Giovanni},
  journal={arXiv preprint arXiv:2503.23832},
  year={2025}
}

@article{ciaperoni2024hadamard,
  title={The {H}adamard decomposition problem},
  author={Ciaperoni, Martino and Gionis, Aristides and Mannila, Heikki},
  journal={Data Mining and Knowledge Discovery},
  volume={38},
  number={4},
  pages={2306--2347},
  year={2024},
  publisher={Springer}
}

@article{oneto2023hadamard,
  title={Hadamard-{H}itchcock decompositions: identifiability and computation},
  author={Oneto, Alessandro and Vannieuwenhoven, Nick},
  journal={arXiv preprint arXiv:2308.06597},
  year={2023}
}

@incollection{friedenberg2017minkowski,
  title={Minkowski sums and {H}adamard products of algebraic varieties},
  author={Friedenberg, Netanel and Oneto, Alessandro and Williams, Robert L},
  booktitle={Combinatorial Algebraic Geometry: Selected Papers From the 2016 Apprenticeship Program},
  pages={133--157},
  year={2017},
  publisher={Springer}
}

@article{hyeon2021fedpara,
  title={Fedpara: Low-rank {H}adamard product for communication-efficient federated learning},
  author={Hyeon-Woo, Nam and Ye-Bin, Moon and Oh, Tae-Hyun},
  journal={arXiv preprint arXiv:2108.06098},
  year={2021}
}

@article{fawzi2015positive,
  title={Positive semidefinite rank},
  author={Fawzi, Hamza and Gouveia, Jo{\~a}o and Parrilo, Pablo A and Robinson, Richard Z and Thomas, Rekha R},
  journal={Mathematical Programming},
  volume={153},
  number={1},
  pages={133--177},
  year={2015},
  publisher={Springer}
}

@article{guglielmi2023rank,
  title={Rank-1 {M}atrix {D}ifferential {E}quations for {S}tructured {E}igenvalue {O}ptimization.},
  author={Guglielmi, Nicola and Lubich, Christian and Sicilia, Stefano},
  journal={SIAM Journal on Numerical Analysis},
  volume={61},
  number={4},
  pages={1737--1762},
  year={2023},
  publisher={SIAM}
}

@article{guglielmi2025low,
  title={A low-rank {ODE} for spectral clustering stability},
  author={Guglielmi, Nicola and Sicilia, Stefano},
  journal={Linear Algebra and its applications},
  volume={721},
  pages={250--276},
  year={2025},
  publisher={Elsevier}
}

@article{guglielmi2024stabilization,
  title={Stabilization of a matrix via a low-rank-adaptive {ODE}},
  author={Guglielmi, Nicola and Sicilia, Stefano},
  journal={BIT Numerical Mathematics},
  volume={64},
  number={4},
  pages={38},
  year={2024},
  publisher={Springer}
}

@book{nesterov2013introductory,
  title={Introductory lectures on convex optimization: A basic course},
  author={Nesterov, Yurii},
  edition = {2nd}, 
  volume={137},
  year={2018},
  publisher={Springer Science \& Business Media}
}

@phdthesis{sicilia2025low,
  title        = {Low-rank properties in structured matrix nearness problems},
  author       = {Sicilia, Stefano},
  year         = 2025,
  month        = {January},
  address      = {},
  note         = {Available at \url{https://iris.gssi.it/handle/20.500.12571/33684?mode=simple}},
  school       = {Gran Sasso Science Institute},
  type         = {Ph{D} {T}hesis}
}

@article{demarinis2025improving,
  title={Improving the robustness of neural {ODE}s with minimal weight perturbation},
  author={De Marinis, Arturo and Guglielmi, Nicola and Sicilia, Stefano and Tudisco, Francesco},
  journal={arXiv preprint arXiv:2501.10740},
  year={2025}
}

@article{guglielmi2017matrix,
  title={Matrix stabilization using differential equations},
  author={Guglielmi, Nicola and Lubich, Christian},
  journal={SIAM Journal on Numerical Analysis},
  volume={55},
  number={6},
  pages={3097--3119},
  year={2017},
  publisher={SIAM}
}

@article{guglielmi2025matrix,
  title={Matrix nearness problems and eigenvalue optimization},
  author={Guglielmi, Nicola and Lubich, Christian},
  journal={arXiv preprint arXiv:2503.14750},
  year={2025}
}

@article{koch2007dynamical,
  title={Dynamical low-rank approximation},
  author={Koch, Othmar and Lubich, Christian},
  journal={SIAM Journal on Matrix Analysis and Applications},
  volume={29},
  number={2},
  pages={434--454},
  year={2007},
  publisher={SIAM}
}

@inproceedings{loconte2022your,
  title={Your {K}nowledge {G}raph {E}mbeddings are {S}ecretly {C}ircuits and {Y}ou {S}hould {T}reat {T}hem as {S}uch},
  author={Loconte, Lorenzo and Di Mauro, Nicola and Peharz, Robert and Vergari, Antonio},
  booktitle={The 5th Workshop on Tractable Probabilistic Modeling},
  year={2022}
}

@article{boumal2014manopt,
  title={Manopt, a {M}atlab toolbox for optimization on manifolds},
  author={Boumal, Nicolas and Mishra, Bamdev and Absil, P-A and Sepulchre, Rodolphe},
  journal={The Journal of Machine Learning Research},
  volume={15},
  number={1},
  pages={1455--1459},
  year={2014},
  publisher={JMLR. org}
}

@article{ang2019accelerating,
  title={Accelerating nonnegative matrix factorization algorithms using extrapolation},
  author={Ang, Andersen Man Shun and Gillis, Nicolas},
  journal={Neural Computation},
  volume={31},
  number={2},
  pages={417--439},
  year={2019},
  publisher={MIT Press One Rogers Street, Cambridge, MA 02142-1209, USA journals-info~…}
}

@article{zhong2005generative,
  title={Generative model-based document clustering: a comparative study},
  author={Zhong, Shi and Ghosh, Joydeep},
  journal={Knowledge and Information Systems},
  volume={8},
  number={3},
  pages={374--384},
  year={2005},
  publisher={Springer}
}

@article{awari2025alternating,
  title={Alternating {D}irection {M}ethod of {M}ultipliers for {N}onlinear {M}atrix {D}ecompositions},
  author={Awari, Atharva and Gillis, Nicolas and Vandaele, Arnaud},
  journal={arXiv preprint arXiv:2512.17473 },
  year={2025}
}

@article{eckart1936approximation,
  title={The approximation of one matrix by another of lower rank},
  author={Eckart, Carl and Young, Gale},
  journal={Psychometrika},
  volume={1},
  number={3},
  pages={211--218},
  year={1936},
  publisher={Springer-Verlag}
}

@article{breslow1983multiplicative,
  title={Multiplicative models and cohort analysis},
  author={Breslow, Norman E and Lubin, Jay H and Marek, P and Langholz, B},
  journal={Journal of the American Statistical Association},
  volume={78},
  number={381},
  pages={1--12},
  year={1983},
  publisher={Taylor \& Francis}
}

@inproceedings{fischer2012introduction,
  title={An introduction to restricted {B}oltzmann machines},
  author={Fischer, Asja and Igel, Christian},
  booktitle={Iberoamerican congress on pattern recognition},
  pages={14--36},
  year={2012},
  organization={Springer}
}

@INPROCEEDINGS{sigmoidharrison,
  author={Nguyen, Harrison and Awari, Atharva and Vandaele, Arnaud and Gillis, Nicolas},
  booktitle={International Workshop on Machine Learning for Signal Processing (MLSP)}, 
  title={Nonlinear {M}atrix {D}ecomposition with the {S}igmoid {F}unction}, 
  year={2025}
  }

@INPROCEEDINGS{seraghitiawarirelu,
  author={Seraghiti, Giovanni and Awari, Atharva and Vandaele, Arnaud and Porcelli, Margherita and Gillis, Nicolas},
  booktitle={International Workshop on Machine Learning for Signal Processing (MLSP)}, 
  title={Accelerated Algorithms For Nonlinear Matrix Decomposition With The {ReLU} Function}, 
  year={2023} 
  }

@article{saulnonlinear,
author = {Saul, Lawrence K.},
title = {A {N}onlinear {M}atrix {D}ecomposition for {M}ining the {Z}eros of {S}parse {D}ata},
journal = {SIAM Journal on Mathematics of Data Science},
volume = {4},
number = {2},
pages = {431-463},
year = {2022},
doi = {10.1137/21M1405769},

URL = { 
    
        https://doi.org/10.1137/21M1405769
    
    

},
eprint = { 
    
        https://doi.org/10.1137/21M1405769
    
    

}
}

@misc{loconte2024turnknowledgegraphembeddings,
      title={How to Turn Your Knowledge Graph Embeddings into Generative Models}, 
      author={Lorenzo Loconte and Nicola Di Mauro and Robert Peharz and Antonio Vergari},
      year={2024},
      eprint={2305.15944},
      archivePrefix={arXiv},
      primaryClass={cs.LG},
      url={https://arxiv.org/abs/2305.15944}, 
}

@INPROCEEDINGS{csfjoakim,
  author={Lefebvre, Joakim and Vandaele, Arnaud and Gillis, Nicolas},
  booktitle={International Workshop on Machine Learning for Signal Processing (MLSP)}, 
  title={Component-{W}ise {S}quared {F}actorization}, 
  year={2024}
  }

@Inbook{onetoalggeom,
author="Friedenberg, Netanel
and Oneto, Alessandro
and Williams, Robert L.",
title="Minkowski {S}ums and {H}adamard {P}roducts of {A}lgebraic {V}arieties",
bookTitle="Combinatorial Algebraic Geometry: Selected Papers From the 2016 Apprenticeship Program",
year="2017",
publisher="Springer New York",
address="New York, NY",
pages="133--157",
isbn="978-1-4939-7486-3",
doi="10.1007/978-1-4939-7486-3_7",
url="https://doi.org/10.1007/978-1-4939-7486-3_7"
}

@article{kileelhad,
title = {Hadamard products of linear spaces},
journal = {Journal of Algebra},
volume = {448},
pages = {595-617},
year = {2016},
issn = {0021-8693},
doi = {https://doi.org/10.1016/j.jalgebra.2015.10.008},
url = {https://www.sciencedirect.com/science/article/pii/S0021869315005323},
author = {Cristiano Bocci and Enrico Carlini and Joe Kileel}
}

@misc{peng2023blockcoordinatedescentsmooth,
      title={Block Coordinate Descent on Smooth Manifolds: Convergence Theory and Twenty-One Examples}, 
      author={Liangzu Peng and René Vidal},
      year={2023},
      eprint={2305.14744},
      archivePrefix={arXiv},
      primaryClass={math.OC},
      url={https://arxiv.org/abs/2305.14744}, 
}

@article{blockroptconv,
  author  = {Yuchen Li and Laura Balzano and Deanna Needell and Hanbaek Lyu},
  title   = {Convergence and complexity of block majorization-minimization for constrained block-{R}iemannian optimization},
  journal = {Journal of Machine Learning Research},
  year    = {2026},
  volume  = {27},
  number  = {42},
  pages   = {1--77},
  url     = {http://jmlr.org/papers/v27/24-0020.html}
}

@article{durrieu2011musically,
  title={A musically motivated mid-level representation for pitch estimation and musical audio source separation},
  author={Durrieu, Jean-Louis and David, Bertrand and Richard, Ga{\"e}l},
  journal={IEEE Journal of Selected Topics in Signal Processing},
  volume={5},
  number={6},
  pages={1180--1191},
  year={2011},
  publisher={IEEE}
}

@inproceedings{durrieu2009main,
  title={Main instrument separation from stereophonic audio signals using a source/filter model},
  author={Durrieu, Jean-Louis and Ozerov, Alexey and F{\'e}votte, C{\'e}dric and Richard, Ga{\"e}l and David, Bertrand},
  booktitle={European Signal Processing Conference},
  year={2009} 
}

@article{durrieu2010source,
  title={Source/filter model for unsupervised main melody extraction from polyphonic audio signals},
  author={Durrieu, Jean-Louis and Richard, Ga{\"e}l and David, Bertrand and F{\'e}votte, C{\'e}dric},
  journal={IEEE Transactions on Audio, Speech and Language Processing}, 
  volume={18},
  number={3},
  pages={564--575},
  year={2010},
  publisher={IEEE}
}
}

\end{document}